\documentclass[10pt]{article}

\usepackage{amsmath,amssymb,amsthm,amsfonts,amstext,amsbsy,amscd}
\usepackage{mathrsfs}
\usepackage{enumerate}

\usepackage{bm}

\usepackage{latexsym}
\usepackage[utf8]{inputenc}
\usepackage{dsfont}
\usepackage{color}

\usepackage[top=4cm, bottom=4.5cm, left=3.5cm, right=3.5cm]{geometry}

\usepackage{hyperref}

\renewcommand{\phi}{\varphi}

\newcommand{\PP}{\mathbb{P}}

\newcommand{\Z}{\mathbb{Z}}

\newcommand{\E}{\mathbb{E}}

\newcommand{\R}{\mathbb{R}}
\newcommand{\N}{\mathbb{N}}

\newcommand{\V}{\mathbb{V}}

\newcommand{\eps}{\varepsilon}

\newcommand{\1}{\mathds{1}}
\newcommand{\peps}{{p}_{\Delta,\eps} }

\newcommand{\nB}{\mathbf{n}(\eps)}
\renewcommand{\hat}{\widehat}
\renewcommand{\tilde}{\widetilde}
\newcommand{\tnB}{\widetilde{\bm{n}}(\eps)}

\newtheorem{theorem}{Theorem}

\newtheorem{rem}{Remark}
\newtheorem{lemma}{Lemma}
\newtheorem{definition}{Definition}

\newtheorem{cor}{Corollary}
\newtheorem{proposition}{Proposition}

\begin{document}

\title{
Spectral-free estimation of Lévy densities in high-frequency regime}

\author{Céline Duval  \footnote{Universit\'e Paris Descartes,  MAP5, UMR CNRS
8145. E-mail: \href{mailto: celine.duval@parisdescartes.fr}{celine.duval@parisdescartes.fr}} \ \ and Ester Mariucci \footnote{Otto von Guericke Universität Magdeburg, Germany.
    E-mail: \href{mailto: mariucci@ovgu.de}{mariucci@ovgu.de}.}}
\date{}

\maketitle

\begin{abstract}
We construct an estimator of the L\'evy density of a pure jump Lévy process, possibly of infinite variation, from the discrete observation of one trajectory at high frequency. The novelty of our procedure is that we directly estimate the Lévy density relying on a pathwise strategy, whereas existing procedures rely on spectral techniques. By taking advantage of a compound Poisson approximation, we circumvent the use of spectral techniques and in particular of the Lévy--Khintchine formula. A linear wavelet estimator is built and its performance is studied in terms of $L_p$ loss functions, $p\geq 1$, over Besov balls. We recover classical nonparametric rates for finite variation Lévy processes and for a large nonparametric class of symmetric infinite variation Lévy processes. We show that the procedure is robust when the estimation set gets close to the critical value 0 and also discuss its robustness to the presence of a Brownian part.\end{abstract}


\section{Introduction}
\subsection{Motivations\label{sec:intro}}

It is now acknowledged that diffusion processes with jumps are good tools for modeling time varying random phenomena whose evolution exhibits sudden changes in value. One of the simplest way to allow for jumps is by considering a Lévy process, that is a continuous time process of the form $X_t=bt + \Sigma W_t+ {\rm Jumps}$, where $W$ is a Brownian motion. The peculiarity of a Lévy process $X$ is that for any $t>0$, the law of $X_t$ is infinitely divisible and the paths of $X$ may have discontinuities. This explains why Lévy processes are a fundamental building block of many stochastic models; many of them have been suggested and extensively studied, for example, in \emph{mathematical finance}; in \emph{physics}, for turbulence, laser cooling and in quantum theory; in \emph{engineering} for networks, queues and dams; in \emph{economics} for continuous time-series models, in \emph{actuarial science} for the calculation of insurance and re-insurance risk (see e.g. \cite{barndorff2012levy,biagini2015electricity,boxma2011levy,carr02,noven2015levy} for reviews and other applications).

 The continuous part of $X$ is characterized by two real parameters $(b,\Sigma^{2})$ and it can be handled easily. The behavior of the jump part is instead described by an infinite-dimensional object, the so-called Lévy measure or, equivalently, by the Lévy density whenever the Lévy measure admits a density with respect to the Lebesgue measure. If the Lévy density $f$ is continuous, $f(x_0)$ determines how frequent jumps of size close to $x_0$ are to occur per unit of time. Thus, to understand the jump behavior of $X$, it is of crucial importance to estimate $f$. 

When dealing with Lévy processes, two approaches are typically used:
\begin{itemize}
\item A spectral approach based on the Lévy--Khintchine formula  
 which relates the cha\-rac\-teristic function of $X_{t}$ to the Lévy density $f$.
 \item A pathwise approach based on the Lévy--Itô decomposition (see \eqref{eq:Xito} below).
\end{itemize}
Techniques employed to address the estimation problem of Lévy densities systematically rely on spectral approaches. They 
 have proven their efficiency both theoretically and numerically. An exception is represented by \cite{Lopez09}, where the properties of a projection estimator of the Lévy density are discussed. In the present work we circumvent the use of spectral techniques in favor of a pathwise strategy. Having a spectral-free procedure paves the way to new techniques for studying richer classes of jump processes for which an equivalent of the Lévy--Khintchine formula is not available. 
 
 One way to proceed is to translate from a probabilistic to a statistical setting Corollary 8.8 in \cite{sato}:
``\emph{Every infinitely divisible distribution is the limit of a sequence of compound Poisson distributions}.'' We are also motivated by the fact that a compound Poisson approximation has been successfully applied to approximate general pure jump Lévy processes, both theoretically and for applications. For example, it is a standard way to simulate trajectories of pure jump Lévy processes (see e.g. Chapter 6, Section 3 in \cite{tankov}). An alternative strategy would consist in taking advantage of the
asymptotic equivalence result in \cite{Mlevy} to construct an estimator of the Lévy density $f$. Yet, the resulting estimator would have the 
strong disadvantage of being randomized and, more fundamentally, it would require the knowledge of $f$ in a neighborhood of the origin.

  In the literature, nonparametric estimation of finite Lévy densities, i.e. Lévy densities of compound Poisson processes, is well understood both from high frequency and low frequency observations (see, among others, \cite{belomestny2018nonparametric,MR2001642,coca2018efficient,coca2018adaptive,Duval,van2007kernel}). 
Building an estimator of $f$ for a Lévy process $X$ with infinite Lévy measure is a more demanding task; for any time interval $[0,t]$, the process $X$ almost certainly jumps infinitely many times. In particular, $f$ is unbounded in any neighborhood of the origin. The techniques used for compound Poisson processes do not generalize immediately. 
Nevertheless, many results on the estimation of -infinite- Lévy densities from discrete data already exist. Spectral techniques enabled to build estimates of functionals of the Lévy density, such as $x f(x)$ or $x^{2} f(x)$ ---which arise naturally when using the Lévy--Khintchine formula--- leading to estimators of $f$ on compact sets away from 0.  A non-exhaustive list of works estimating $f$ for $L_2$ and $L_\infty$ loss functions includes \cite{MR2816339,MR2565560,figueroa2011sieve,LH,gugushvili2009nonparametric,Shota12,kappus2015nonparametric,NR09,Trabs15}; a review is also available in the textbook \cite{BCGM}.

\subsection{Notations and definitions\label{sec:not}}

Before detailing the results, we introduce some necessary notations and definitions.
\paragraph*{Lévy--Itô decomposition} 

It is well known that any Lévy process $X$ is a càdlàg process that can be written as
\begin{align}\label{eq:Xgen}
X_t&=bt+\Sigma W_t+  \lim_{\eta\to 0}\bigg(\sum_{s\leq t}\Delta X_s\1_{|\Delta X_s|\leq 1}-t\int_{|x|\leq 1}x\nu(dx)\bigg)+ \sum_{0<s\leq t}\Delta X_s \1_{|\Delta X_s|> 1}\nonumber \\
&=:bt+\Sigma W_t+M_t(1)+Z_t(1),
\end{align}
where $b\in \R$, $\Sigma\in \R_{\geq 0}$, $\nu$ is a Borel measure on $\R$ such that $\nu(\{0\})=0$ and $ \int_{\R}(y^2\wedge 1)\nu(dy)<\infty$, $W=(W_t)_{t\geq 0}$ is a standard Brownian motion, $\Delta X_r=X_r-\lim_{s\uparrow r} X_s$ and the processes $W$, $M$ and $Z$ are independent.

In this paper we focus on the class of pure jump Lévy processes of the form
\begin{equation}\label{eq:defX} X_t:=\begin{cases}
             \sum_{0<s\leq t}\Delta X_s\quad &\text{if}\quad \int_{|x|\leq 1}|x|\nu(dx)<\infty,\\
               M_t(1)+Z_t(1)\quad &\text{if}\quad \int_{|x|\leq 1}|x|\nu(dx)=\infty.
             \end{cases}   
\end{equation}

For all $0< \varepsilon \le 1$, the Lévy--Itô decomposition allows to write a Lévy process $X$ as in \eqref{eq:defX} as the sum of two independent Lévy processes, the first one (resp. the second one) having jumps smaller (resp. larger) in absolute value than $\eps$. Namely
\begin{align}\label{eq:Xito}
 X_t&=tb_\nu(\varepsilon)+\lim_{\eta\to 0}\bigg(\sum_{s\leq t}\Delta X_s\1_{(\eta,\varepsilon]}(|\Delta X_s|)-t\int_{\eta<|x|\leq \varepsilon}x\nu(dx)\bigg) + \sum_{i=1}^{N_t(\varepsilon)}Y_i(\varepsilon)\nonumber \\&=: tb_\nu(\varepsilon)+M_t(\varepsilon)+Z_t(\varepsilon),
\end{align}
where the drift $b_\nu(\varepsilon)$ is defined as 
\begin{equation}\label{eq:bnu}b_\nu(\varepsilon):=\begin{cases}
             \int_{|x|\leq \varepsilon} x\nu(dx)\quad &\text{if}\quad \int_{|x|\leq 1}|x|\nu(dx)<\infty,\\
                -\int_{\varepsilon\leq |x|\leq 1} x\nu(dx)\quad &\text{if}\quad \int_{|x|\leq 1}|x|\nu(dx)=\infty ,
             \end{cases}   
\end{equation}
$M(\varepsilon)=(M_t(\varepsilon))_{t\geq 0}$ and $Z(\varepsilon)=(Z_t(\varepsilon))_{t\geq 0}$ are two independent Lévy processes. The process  $M(\varepsilon)$ is a centered martingale consisting of the sum of the \emph{small jumps} i.e. the jumps of size smaller than $\varepsilon$. The process 
 $Z(\varepsilon)$ instead, is a compound Poisson process defined as follows: $N(\varepsilon)=(N_t(\varepsilon))_{t\geq 0}$ is a Poisson process of intensity $\lambda_\varepsilon:=\int_{|x|>\varepsilon}\nu(dx)$ and $(Y_i(\varepsilon))_{i\geq 1}$ are i.i.d. random variables independent of $N(\varepsilon)$ such that $\PP(Y_1(\varepsilon)\in A)={\nu(A)}/{\lambda_\varepsilon}$, for all $A\in\mathscr B(\R\setminus (-\varepsilon,\varepsilon))$.

\paragraph*{Compound Poisson approximation}

Denote by $f$ (resp. $f_\varepsilon$) the \emph{Lévy density} of $X$ (resp. $Z(\varepsilon)$), i.e. $f(x)=\frac{\nu(dx)}{dx}$ (resp. $f_\varepsilon(x)=\frac{\1_{|x|>\varepsilon}\nu(dx)}{dx}$). Let $h_\varepsilon$ be the density, with respect to the Lebesgue measure, of the random variables $(Y_i(\varepsilon))_{i\geq 0}$, i.e. $h_{\eps}=f_{\eps}/\lambda_{\eps}$. 
We are interested in estimating $f$ in any set of the form $A(\varepsilon):=(-\overline{A},-\varepsilon]\cup [\varepsilon,\overline{A})$ for all $\eps>0$, where $\overline{A}\in (\eps,\infty]$. The latter condition is technical; if $X$ is a compound Poisson process we may choose $\overline{A}:=+\infty$, otherwise we work under the 
simplifying assumption that $A(\eps)$ is a 
bounded set. Observe that, for any $\eps>0$,
\begin{equation}\label{eq:f}
 f(x) \1_{A(\eps)}(x)=f_\eps(x)\1_{A(\eps)}(x)= \lambda_\varepsilon h_\varepsilon(x) \1_{A(\eps)}(x),\  \forall x\in\R.
\end{equation} 
Therefore, estimating $f$ in $A(\eps)$ from the increments of $X$ is equivalent to estimating the Lévy density of the compound Poisson part of $X$, namely $Z(\eps)$, from the increments of $X$.
When $\nu(\R)<\infty$, we may take $\eps=0$ and Equation \eqref{eq:Xito} reduces to a compound Poisson process
with intensity $\lambda =\nu(\R)$ and jump density $h=f/\lambda$.

\paragraph*{A nonparametric class of Lévy densities}   Assume  that the Lévy measure $\nu$ is absolutely continuous with respect to the Lebesgue measure and denote by $f$ the \emph{Lévy density} of $X$. We pay special attention to the following nonparametric class of Lévy densities. Consider $\alpha\in(0,2)$ and $M$ a positive constant, define the class of functions
\begin{align}\label{eq:Lalpha}
\mathscr L_{M,\alpha}&:=\bigg\{f : f(x)\leq \frac{M}{|x|^{1+\alpha}},\quad  \forall |x|\leq 2\bigg\}.
\end{align}
A Lévy density $f$ belongs to the class  $\mathscr L_{M,\alpha}$, $M>0$, $\alpha\in(0,2)$, if $\sup_{x\in[-2,2]}f(x)|x|^{1+\alpha}\le M$. In particular $\mathscr L_{M,\alpha}$ contains any $\tilde\alpha$-stable Lévy density such that $\tilde\alpha\le\alpha$. Any finite variation Lévy process is in the class $\mathscr	L_{M,1},$ for some positive $M$. 

In Theorem \ref{thm:main}, we provide upper bounds for general Lévy densities. To derive explicit rates of convergence we will examine, in particular, the cases where $f$ belongs to $\mathscr L_{M,\alpha}$ relying on the results of \cite{DMnote} (see Appendix \ref{app:note}).

\paragraph*{Observation setting and loss function}
Suppose we observe $X$ on $[0,T]$ at the sampling rate $\Delta>0$. Without loss of generality, set $T:=n\Delta$ with $n\in \N$, and define
\begin{equation}\label{eq:xi}
\mathbf{X}_{n,\Delta}:=
(X_\Delta, X_{2\Delta}-X_{\Delta},\dots, X_{n\Delta}-X_{(n-1)\Delta}).
\end{equation}
We consider the high frequency setting where $\Delta\to 0$ and $T=n\Delta\to\infty$ as $n\to\infty$. The assumption $n\Delta\to\infty$ is necessary to construct a consistent estimator of $f$. 
The difference with the works listed in Section \ref{sec:intro} is that we build a spectral-free estimator of $f$, without smoothing treatment at the origin, and study the following $L_{p}$ risk. Define the class $L_{p,\varepsilon}=\big\{g:\|g\|_{L_p,\varepsilon}:=\Big(\int_{A(\varepsilon)} |g(x)|^pdx\Big)^{{1}/{p}}<\infty\big\},$ where $ p\in[1,\infty)$ and $A(\varepsilon)$ is the estimation set defined above.
Define the loss function
$$
\ell_{p,\varepsilon}\big(\widehat f,f\big):=\big(\E\big[\big\|\widehat f-f\big\|_{L_{p,\varepsilon}}^{p}\big]\big)^{1/p}
=\bigg(\E\bigg[\int_{A(\varepsilon)}|\hat f(x)-f(x)|^pdx\bigg]\bigg)^{1/p},\quad p\in[1,\infty).
$$

  Finally, denote by $P_\Delta$ the distribution of the random variable $X_\Delta$ and by $P_n$ the law of the random vector $\mathbf{X}_{n,\Delta}$ defined in \eqref{eq:xi}. Since $X$ is a Lévy process, its increments are i.i.d., hence 
\begin{equation*}
 P_n=\bigotimes_{i=1}^n P_{i,\Delta}=P_{\Delta}^{\otimes n},\quad \textnormal{where}\quad P_{i,\Delta}=\mathscr L(X_{i\Delta}-X_{(i-1)\Delta}).
\end{equation*}
 In the following, whenever confusion may arise, the reference probability in expectations is explicitly stated, for example, writing $\E_{P_n}$.

\subsection{Estimation strategy and results} For any fixed $0<\varepsilon\le 1$ (when $\nu(\R)<\infty$ the choice $\varepsilon=0$ is allowed), taking advantage of Equation \eqref{eq:f}, we build an estimator of $f$ on the set $A(\varepsilon)$ by constructing estimators for $\lambda_\varepsilon$ and $h_\varepsilon$ separately.  For that we consider the increments of \eqref{eq:xi} larger than $\eps$ in absolute value. Define the dataset \begin{align}\label{eq:data}\mathbf{D}_{n,\varepsilon}:=\big\{X_{i\Delta}-X_{(i-1)\Delta},\ i\in \mathscr I_\varepsilon\big\},\end{align}
 where $ \mathscr I_\varepsilon$ is the subset of indices such that
$ \mathscr I_\varepsilon:=\big\{i=1,\dots,n:|X_{(i-1)\Delta}-X_{i\Delta}|> \varepsilon\big\}$. Its random cardinality is denoted by
\begin{equation}\label{eq:neps}
 \mathbf{n}(\varepsilon):=\sum_{i=1}^n\1_{\R\setminus[-\varepsilon,\varepsilon]}(|X_{i\Delta}-X_{(i-1)\Delta}|).
\end{equation}
Our estimation strategy is the following.
\begin{enumerate}
\item We build an estimator of $\lambda_{\eps} $ using the following Lemma:
\begin{lemma}\label{lemma:momenti} Let $X$ be a Lévy process with Lévy measure $\nu$ absolutely continuous with respect to the Lebesgue measure. Then, for all $\eps\in(0,1]$
\begin{equation*}
\lim_{t\to0}\frac{ \PP(|X_t|\geq \eps)}{t}=\int_{\R\setminus[-\eps,\eps]} \nu(dy).
\end{equation*}
\end{lemma} In particular, Lemma \ref{lemma:momenti} implies 
$
\lim_{t\to0}t^{-1}{\PP(M_t(\eps)\geq\eps)}=0
$ and \begin{equation}\label{eq:l}
\lambda_{\eps}=\lim_{t\rightarrow 0}\frac{1}{t}\PP( |X_{t}|> \eps), \quad \forall \  \eps\in(0,1].
\end{equation}
Lemma \ref{lemma:momenti} is a modification of Lemma 6 in Rüschendorf and Woerner \cite{Rusch2002}. Their proof relies on spectral arguments, we provide a spectral free version of the proof in the Appendix.
 \item From the observations $\mathbf{D}_{n,\eps}$ in \eqref{eq:data} we build a wavelet estimator $\hat h_{n,\eps}$ of  $h_{\eps}$ using that, for $\Delta$ small, the random variables $(X_{i\Delta}-X_{(i-1)\Delta})_{i\in\mathscr I_\varepsilon}$ are i.i.d. with a density close to $h_{\eps}$ (see Lemma \ref{lem:hpeps} below). 
\item Finally, we estimate $f$ on $A(\eps)$ following \eqref{eq:f} by \begin{equation}\label{eqhatf}
 \widehat f_{n,\varepsilon}(x):=\widehat \lambda_{n,\eps}\widehat h_{n,\varepsilon}(x)\1_{A(\varepsilon)}(x),\quad \forall x\in A(\eps).
\end{equation}
\end{enumerate}

The presence of the small jumps makes the estimation of  $h_{\eps}$ and $\lambda_{\eps}$ from observed increments larger than $\eps$ delicate. Indeed, if $i_{0}$ is such that $|X_{i_{0}\Delta}-X_{(i_{0}-1)\Delta}|>\varepsilon$,  it is not automatically true that there exists $s\in((i_{0}-1)\Delta,i_{0}\Delta]$ such that $|\Delta X_{s}|>\varepsilon$ (or any other fixed 
positive number).

In Section \ref{sec:main} we establish upper bounds for the $L_p$ risks of the estimators $\widehat \lambda_{n,\eps}$ and $\widehat h_{n,\varepsilon}$ (see Theorem \ref{teolambda},  Proposition \ref{teo:hB} and Corollary \ref{cor:hB}). The main difficulty in their study lies in the presence of the small jumps that play  a role in both cases.  We stress that when $\eps$ is fixed, the quantity $\lambda_{\eps}$ is bounded. But as we generalize the results to the case $\eps\to 0$ (see Section \ref{sec:eps0}) this is no longer true: whenever $\nu(\R)=\infty$, it holds $\lambda_{\eps}\to\infty$ as $\eps\to0$. Therefore, in all the results of the paper we always keep explicit the dependency in $\lambda_{\eps}$.

Our main results on the estimator \eqref{eqhatf} of the Lévy density $f$ are given in Section \ref{sec:mainthm}. Theorem \ref{thm:main} provides a general upper bound that tends to 0, regardless of the rate at  which $\Delta$ tends to 0. Interestingly, the estimation strategy leads to an upper bound where terms depending on the behavior of the small jumps appear. These terms make it difficult to derive an explicit rate of convergence without additional assumptions on the Lévy density. Therefore, in Theorems \ref{thm:mainF} and \ref{thm:mainF2} we consider additional assumptions, satisfied in particular by the classes introduced in \eqref{eq:Lalpha}, and derive explicit rates of convergence. Similarly to what happens when using spectral procedures, Theorem \ref{thm:mainF} ensures that we recover the classical rates for finite variation Lévy densities when $n\Delta^{2}\le 1$ (see e.g. \cite{MR2565560,MR2816339}). Furthermore, we show that this rate is attained for a large nonparametric class of symmetric infinite variation Lévy processes. Theorem \ref{thm:mainF2} generalizes the rates of Theorem \ref{thm:mainF}, in particular to slower regimes for $\Delta$. Finally, Theorem \ref{thm:mainF3} shows that our procedure is robust to the case $\eps\to 0$. 
For the sake of clarity, the results of Sections \ref{sec:main} and \ref{sec:mainthm} are stated for purely discontinuous Lévy processes, even though our procedure generalizes to the presence of a Gaussian part as detailed in Section \ref{sec:brown}. 

 Finally, Section \ref{sec:prf} contains the proofs of the main results while Appendix \ref{sec:app} collects the proofs of the auxiliary results and in Appendix \ref{app:note} technical results of \cite{DMnote} are partially reproduced and used to establish Theorems \ref{thm:mainF}, \ref{thm:mainF2} and \ref{thm:mainF3}.

\section{Preliminary estimators\label{sec:main}}

\subsection{Statistical properties of $\hat\lambda_{n,\eps}$}\label{sec:lbd}

First, we define an estimator of the intensity of the Poisson process $Z(\eps)$ in terms of $\nB$, the number of jumps that exceed $\eps$ following  \eqref{eq:l}.
\begin{definition}\label{def:estlbd}
Let $\widehat \lambda_{n,\eps}$ be the estimator of $\lambda_{\eps}=\int_{|x|>\eps}\nu(dx)$ defined by 
\begin{align}\label{eq:lambdahat}
\widehat \lambda_{n,\eps}:=\frac{\mathbf{n}(\eps)}{n\Delta},
\end{align}
where $\mathbf{n}(\varepsilon)$ is defined as in \eqref{eq:neps}.\end{definition}
Observe that  $ \hat \lambda_{n,0}\asymp N_{n\Delta}/(n\Delta)$ as  $\Delta\to 0$,  which is the maximum likelihood estimator of $\lambda_{0}$ in the experiment $\{\mathcal{P}(\lambda_{0}\Delta),\ \lambda_{0}\in(0,\infty)\}$.
 In \cite{NR12} and \cite{NRST16},  estimators of the cumulative distribution function of the Lévy measure, which is closely related to $\lambda_\varepsilon$, are build and Donsker theorems are derived. In \cite{NRST16} a direct approach similar to \eqref{eq:lambdahat} is considered, the performances are investigated in $L_{\infty}(V)$, for a domain $V$ bounded away from 0.
 We establish the following $L_{p}$ bound for $\widehat\lambda_{n,\eps}$.

\begin{theorem}\label{teolambda} Let $X$ be a Lévy process as in \eqref{eq:defX} and let
 $n\geq 1$, $\Delta>0$, $\eps\in(0,1]$ be such that $n\PP(|X_\Delta|>\varepsilon)\geq 1$. Let $\widehat \lambda_{n,\eps}$ be the estimator of $\lambda_\eps$ defined in \eqref{eq:lambdahat}. 
Then, there exists a constant $C$, depending only on $p$, such that
 \begin{align*}
&\E_{P_n}\big[|\widehat \lambda_{n,\eps}-\lambda_\varepsilon|^p\big]\leq \Big|\lambda_\eps-\frac{\PP(|X_\Delta|>\varepsilon)}{\Delta}\Big|^p
+C \Big(\frac{\PP(|X_\Delta|>\varepsilon)}{n\Delta^2}\Big)^{\frac{p}{2}}, \quad \forall p\in [1,\infty).
\end{align*}
\end{theorem}
In general, the quantity ${\PP(|X_\Delta|>\eps)}/{\Delta}$ is not easy to handle. By Lemma \ref{lemma:momenti}, it holds
$\lim_{\Delta\to 0}\bigg|\lambda_\eps-\frac{\PP(|X_\Delta|>\eps)}{\Delta}\bigg|=0,$
but the rate of convergence is not known in general. Nevertheless, in many cases of interest, it holds that
 $|\lambda_\eps-\frac{\PP(|X_\Delta|>\eps)}{\Delta}|=O(\Delta \lambda_\eps^2),$ as $ \Delta\lambda_\eps\to 0.$
This motivates Assumption \eqref{eq:assH1} below  that leads to Corollary \ref{cor:teolambda}.

\smallskip
 
\textbf{Assumption $H_1(\delta,c)$}: $X$ is a Lévy process with a Lévy measure $\nu$ such that
\begin{equation}
\big|\PP(|X_t|>\eps)-t\lambda_\eps\big|\le ct^2,\quad \forall 0<t\leq \delta,\ 0<\eps\le 1,\tag{$H_1(\delta,c)$}\label{eq:assH1}
\end{equation}where $c$ is a constant that does not depend on $t$.
\begin{cor}\label{cor:teolambda} Let $X$ be as in \eqref{eq:defX} and such that \eqref{eq:assH1} is satisfied for some $\delta>0$ and $c>0$. Let $0<\Delta\leq\delta$, $\eps\in(0,1]$ and $n\ge 1$; the estimator \eqref{eq:lambdahat} of $\lambda_{\eps}$ satisfies for all $p\in [1,\infty)$
 \begin{align*}
&\E_{P_n}\big[|\widehat \lambda_{n,\eps}-\lambda_{\varepsilon}|^p\big]\le C\bigg(\Big(\frac{\lambda_{\eps}+\Delta}{n\Delta}\Big)^{\frac{p}{2}}+\Delta^p\bigg),\end{align*} where $C$ is a positive constant depending on $p$ and $c$.
\end{cor}
Assumption \eqref{eq:assH1} requires a non-asymptotic control on the cumulative distribution function of the Lévy process for small times. Asymptotic expansions  have been established in the literature such as Theorem 3.2 in \cite{LH}, where a control  of $\PP(X_t>y)$ is given for $y$ bounded away from the origin. However, no indication on how small $t$ should be nor  on how large $y$ should be is given. 

It is possible to establish that \eqref{eq:assH1} holds true on the nonparametric class $\mathscr L_{M,\alpha}$	(see \eqref{eq:Lalpha}). This is a consequence of \cite{DMnote} whose main results are reproduced in Appendix \ref{app:note}. Theorem \ref{thm:Fd1} ensures that for $\alpha\in(0,1)$ and $M>0$, \eqref{eq:assH1} is satisfied for $\delta=\frac{(1-\alpha)\eps^{\alpha}}{M4^{1+\alpha}}$ and $c$ depending on $\alpha,\ \eps$ and $M$. Moreover, if $X$ is symmetric and its Lévy density is $M\eps^{-(2+\alpha)}$-Lipschitz on the interval $(3/4\eps,5/4\eps)$, Theorem \ref{thm:Fd2} ensures that for $\alpha\in[1,2)$ and $M>0$, \eqref{eq:assH1} is satisfied for $\delta=\frac{(2-\alpha)\eps^{\alpha}}{M2^{1+\alpha}}$ and $c$ depending on $\alpha,\ \eps$ and $M$. Note that any L\'evy density of the form $f(x)=L(x)/x^{1+\alpha}$ for all $x\in[-2,2]\setminus \{0\}$ where $L$ is a bounded differentiable function with bounded derivative is $M\eps^{-(2+\alpha)}$-Lipschitz on the interval $(3/4\eps,5/4\eps)$. In particular it is satisfied if $X$ is an $\alpha$-stable process.
 Furthermore, in Theorems \ref{thm:Fd1} and \ref{thm:Fd2}, the dependency in $\eps$ of the constant $c$ is explicit, which will allow to control the $L_{p}$ loss between  $\hat \lambda_{n,\eps}$ and $\lambda_{\eps}$ in the asymptotic $\eps\to0$.

\subsection{Statistical properties of $\widehat h_{n,\eps}$\label{sec:h}}

\subsubsection{Construction of $\widehat h_{n,\eps}$}
To recover the jump density $h_\varepsilon$, we exploit the high frequency setting. For $\Delta$ small, it holds $h_{\eps}\approx\mathcal{L}(X_{\Delta}\big||X_{\Delta}|>\eps)$. Focusing on the increments larger than $\eps$ in absolute value, we estimate the density $h_\varepsilon$ using a linear wavelet
density estimator and study its performances uniformly over Besov balls (see Kerkyacharian and Picard \cite{kerkyacharian1992density} or H\"ardle et al. \cite{hardle2012wavelets}). We state the result and assumptions in terms of the quantity of interest $f$.

\paragraph{Preliminary on Besov spaces}

Let $(\Phi,\Psi)$ be a pair of scaling function and mother wavelet which are compactly supported, of class $C^{r}$ and generate a regular wavelet basis
adapted to the estimation set $A(\eps)$ (e.g. Daubechie's wavelet). Moreover suppose that $\{\Phi(x-k),k\in \Z\}$ is an orthonormal family of $L_{2}(\R)$. For all $f\in L_{p,\eps}$  we write for $j_{0}\in\N$
\begin{equation}
f(x)=\sum_{k\in\Lambda_{j_{0}}}\alpha_{{j_{0}}k}(f)\Phi_{{j_{0}}k}(x)+\sum_{j\geq {j_{0}}}\sum_{k\in\Lambda_j}\beta_{jk}(f)\Psi_{jk}(x),\quad\forall x\in A(\eps),\nonumber
\end{equation}
where $\Phi_{j_{0}k}(x)=2^{\frac{j_{0}}{2}}\Phi(2^{j_{0}}x-k)$, $\Psi_{jk}(x)=2^{\frac{j}{2}}\Psi(2^jx-k)$ and the coefficients are
\begin{equation*}
 \alpha_{j_{0}k}(f)=\int_{A(\eps)}\Phi_{j_{0}k}(x)f(x)dx\quad\text{and}\quad \beta_{jk}(f)=\int_{A(\eps)}
 \Psi_{jk}(x)f(x)dx.
\end{equation*} As we consider compactly supported wavelets, for every $j\geq j_{0}$, the set $\Lambda_j$  incorporates boun\-da\-ry terms that we choose not to distinguish in notation 
for simplicity. In the sequel we apply this decomposition to $h_{\eps}$. This is justified because $f_{\eps}\in L_{p,\eps}$ implies $h_{\eps}\in L_{p,\eps}$ and the coefficients of its decomposition are $\alpha_{j_{0}k}(h_{\eps})=\alpha_{j_{0}k}(f)/\lambda_{\eps}$ and $\beta_{j_{0}k}(h_{\eps})=\beta_{j_{0}k}(f)/\lambda_{\eps}$. The latter can be interpreted as the expectations of $\Phi_{j_{0}k}(U)$ and $\Psi_{jk}(U)$ where $U$ is a random variable with density $h_{\eps}$ with respect to the Lebesgue measure.

 We define Besov spaces in terms of wavelet coefficients as follows. 
For $r>s>0$, $p\in[1,\infty)$ and $1\leq q\leq\infty$ a function $f$ belongs
to the Besov space $B_{p,q}^{s}(A(\eps))$ if the norm
{\small\begin{align}
 \|f\|_{B_{p,q}^{s}(A(\eps))}:=\bigg(\sum_{k\in\Lambda_{j_{0}}}|\alpha_{j_{0}k}(f)|^p\bigg)^{\frac{1}{p}}+\bigg[\sum_{j\geq j_{0}}\bigg(2^{j(s+1/2-1/p)}\Big(\sum_{k\in\Lambda_j}|
 \beta_{jk}(f)|^p\Big)^{\frac{1}{p}}\bigg)^{q}\bigg]^{\frac1q}\label{eq:Besovnorm}
\end{align}}
is finite, with the usual modification if $q=\infty$.
We consider Lévy densities $f$ with respect to the Lebesgue measure, whose restriction to the set $A(\eps)$ lies into a Besov ball:
\begin{equation}\label{eq:Fbesov}
 \mathscr F({s,p,q,\mathfrak{M}_\eps,A(\eps)})=\big\{f\in L_{p,\eps}\ :\ \|f\|_{B_{p,q}^{s}(A(\eps))}\leq \mathfrak M_{\eps}\big\},
\end{equation} where $\mathfrak{M}_\eps:=\mathfrak{M}\lambda_\eps$, for a fixed constant $\mathfrak{M}$.
Note that the regularity assumption is imposed on $f_{|A(\eps)}$ viewed as an $L_{p,\eps}$ function. Therefore the dependency in $A(\eps)$ lies in $\mathfrak M_{\eps}$. Also, the parameter $p$ measuring the loss of our estimator is the same as the one measuring the Besov regularity of the function, this is discussed in Section \ref{sec:UBR}. Lemma \ref{lem:hepsBesov} below follows immediately from the definitions of $h_{\eps}$ and the Besov norm \eqref{eq:Besovnorm}.

\begin{lemma}\label{lem:hepsBesov}
For all $0<\eps\le 1$, let $f$ be in $ \mathscr F({s,p,q,\mathfrak{M}_\eps,A(\eps)})$. Then, $h_{\eps}=\frac{f_{\eps}}{\lambda_{\eps}}$ belongs to the class $\mathscr F\big(s,p,q,{\mathfrak	M},A(\eps)\big)$.
\end{lemma}

\paragraph{Construction of $\hat h_{n,\eps}$}
To estimate the jump density $h_{\eps}$, we only have access to the indirect observations $\{X_{i\Delta}-X_{(i-1)\Delta},\ i\in \mathscr{I}_{\eps}\}$, where
for each $i\in \mathscr I_{\eps}$, it holds
 $$X_{i\Delta}-X_{(i-1)\Delta}=M_{i\Delta}(\eps)-M_{(i-1)\Delta}(\eps)+\Delta b_{\nu}(\eps) +Z_{i\Delta}(\eps)-Z_{(i-1)\Delta}(\eps).$$
 The problem is twofold. First, there is a deconvolution problem as the information on $h_{\eps}$ is contained in the observations $\{Z_{i\Delta}(\eps)-Z_{(i-1)\Delta}(\eps), i\in \mathscr I_{\eps}\}$. The distribution of the noise $M_{\Delta}(\eps)+\Delta b_{\nu}(\eps) $ is unknown, but it is small as its variance $\V(M_{\Delta}(\eps))=\Delta\int_{|x|\leq \eps}x^{2}\nu(dx)\to0$ as $\Delta\to0$. Then, we neglect this noise:
\begin{align}
\label{eq:M0}X_{i\Delta}-X_{(i-1)\Delta}\approx Z_{i\Delta}(\eps)-Z_{(i-1)\Delta}(\eps),\quad \forall	i\in \mathscr I_{\eps}.\end{align}
 Second, overlooking that it is possible that for some $i_{0}\in\mathscr I_{\eps}$, $|X_{i_{0}\Delta}-X_{(i_{0}-1)\Delta}|>\eps$ and $Z_{i_{0}\Delta}-Z_{(i_{0}-1)\Delta}=0,$ yet the common density of $Z_{i\Delta}-Z_{(i-1)\Delta}|Z_{i\Delta}-Z_{(i-1)\Delta}\ne0$ is not $h_{\eps}$ but it is given by
\begin{align}
\label{eq:peps}\peps &=\sum_{k=1}^{\infty}\PP(N_{\Delta}(\eps)=k|N_{\Delta}(\eps)\ne 0)h_{\eps}^{\star k}=\sum_{k=1}^{\infty}\frac{(\lambda_{\eps}\Delta)^{k}}{k!(e^{\lambda_{\eps}\Delta}-1)}h_{\eps}^{\star k},
\end{align}where $\star$ denotes the convolution product. Again, in the asymptotic $\Delta\rightarrow 0$, we neglect the possibility that more than one jump of $N(\eps)$ occurred in an interval of length $\Delta$. 
\begin{lemma}\label{lem:hpeps}
For all $p\geq 1$, $\eps\in(0,1]$ and $\Delta>0$, it holds 
$\big\|\peps -h_{\eps}\big\|_{L_{p,\eps}}\hspace{-0.3cm}\leq 2\Delta e^{\lambda_\eps\Delta}\|f\|_{L_{p,\eps}}.$
\end{lemma} \noindent Define the estimator based on the chain of approximations $h_{\eps}\approx \peps \approx \mathcal{L}(X_{\Delta}||X_{\Delta}|>\eps)$ 
\begin{align}
\label{eq hat heps}
\hat h_{n,\eps}(x)=\sum_{k\in \Lambda_{J}}\widehat \alpha_{J,k}\Phi_{Jk}(x),\quad x\in A(\eps),
\end{align} where $J$ is an integer to be chosen and
\begin{align}
 \widehat \alpha_{J,k}&:=\frac{\mathds{1}_{\nB\geq 1}}{\nB}\sum_{i\in \mathscr{I}_{\eps}}\Phi_{Jk}( X_{i\Delta}-X_{(i-1)\Delta}).\nonumber
 \end{align} If $\nB=0$, the estimator $\hat h_{n,\eps}$ is 0, which occurs with probability $(1-\PP(|X_{\Delta}|>\eps))^{n}\le e^{-n\PP(|X_{\Delta}|>\eps)}$.  We work with a linear estimator even if linear estimators are not always minimax for general Besov spaces $B^{s}_{\pi,q}$, $1\leq \pi,q\leq \infty$ ($\pi\ne p$). Indeed, to evaluate the loss caused by neglecting the small jumps $M_{\Delta}(\eps)$ (see \eqref{eq:M0}), we make an approximation at order 1 of our estimator $\hat h_{n,\eps}$. We thus require our estimator to depend smoothly on the observations, which is not the case for usual thresholding methods. Finally, we recall that on the class $\mathscr F(s,p,q,\mathfrak M_{\eps},A(\eps))$ this estimator is optimal in the context of density estimation from direct i.i.d. observations (see  \cite{kerkyacharian1992density}, Theorem 3).

\subsubsection{Upper bound results\label{sec:UBR}}

Adapting the results of \cite{kerkyacharian1992density}, we derive a conditional upper bound for the estimation of $h_\eps$ when the Lévy measure is infinite. Recall that $A(\eps)=(-\overline{A},-\eps]\cup[\eps,\overline{A})$ with $\overline{A}\in(\eps,\infty{)}$.

 \begin{proposition}\label{teo:hB}
Suppose $\mathscr{I}_\eps\ne\emptyset$, fix $0<\varepsilon\le 1$ and that $f$ belongs to the class $ \mathscr F({s,p,q,\mathfrak{M}_\eps,A(\eps)})$ defined in \eqref{eq:Fbesov}, for some $1\leq q\leq \infty$, $ 1\leq p<\infty$,  and $\eps<\overline{A}<\infty$. If $1\leq p<2$ suppose that $h_{\eps}(x)\leq w_{\eps}(x)$, $\forall x\in A(\eps)$ for some symmetric function $w_{\eps}\in L_{p/2}$.
 Let $r>s>\frac{1}{p}$ and let $\widehat h_{n,\varepsilon}$ be the wavelet estimator of 
 $h_{\eps} $ on $A(\eps)$, defined in \eqref{eq hat heps}. 
 Let $v_{\Delta}(\eps):=\PP(|M_{\Delta}(\varepsilon)+\Delta b_\nu(\varepsilon) |>\eps)$, $F_{\Delta}(\eps):=\PP(|X_{\Delta}|>\eps)$, $\sigma^2(\varepsilon):=\int_{|x|\leq \varepsilon}x^2\nu(dx)$ and $\mu_p(\eps):=\int_{|x|\leq \varepsilon}|x|^p\nu(dx)$. 
For any $\Delta>0$ such that $\tfrac{v_{\Delta}(\eps)}{F_{\Delta}(\eps)}\leq\tfrac13$,  the following inequality holds. For all $J\in \N$ and for all finite $p\geq 1$, 
\begin{align*}
\E\big[\|\hat h_{n,\eps}(\{X_{i\Delta}&-X_{(i-1)\Delta}\}_{i\in \mathscr I_{\eps}})-h_{\eps}\|_{L_{p,\eps}}
^{p}|\mathscr I_{\eps}\big]\\ &\leq C\bigg\{{2^{2Jp}}\Big[\Big(\frac{v_{\Delta}(\eps){e^{-\lambda_\eps\Delta}}}{\nB F_{\Delta}(\eps)}\Big)^{p/2}+\Big(\frac{v_{\Delta}(\eps){e^{-\lambda_\eps\Delta}}}{F_{\Delta}(\eps)}\Big)^{p}\Big]\\
 & \hspace{-1cm}+2^{-Jsp}+2^{Jp/2}\ell_{p,\eps}^{p/2}\nB^{-p/2}+2^{J(p-1)}\nB^{(1-p)}\mathbf{1}_{p\ge 2}+(e^{\lambda_\eps\Delta}\Delta\|f\|_{L_{p,\eps}})^{p}
 \\
 & \hspace{-2cm}+{2^{J(5\frac p2-1)}}\Big[\nB^{1-p}(\Delta \mu_p(\eps)+(\Delta\sigma^2(\eps))^{\frac p2})+\nB^{-\frac p2}\big(\sigma^{2}(\eps)\Delta\big)^{\frac p2}+(b_\nu(\varepsilon)\Delta)^p\Big]\bigg\},\end{align*} 
where  $\nB$ denotes the cardinality of $\mathscr I_{\eps}$,  $\ell_{p,\eps}:=\|h_{\eps}\|_{L_{p/2,\eps}}\mathbf{1}_{p\ge 2}+\|w_{\eps}\|_{L_{p/2,\eps}\mathbf{1}_{1\le p<2}}$ and $C$ only depends on $s$, $p$, $\|\Phi\|_{\infty}$, $\|\Phi'\|_{\infty}$, $\|\Phi\|_{p}$ and $\mathfrak M$.

\end{proposition}
\paragraph{Comments}
If $\mathscr{I}_\eps=\emptyset$, then $\hat h_{n,\eps}=0$ and we get $\E\big[\|\hat h_{n,\eps}-h_{\eps}\|_{L_{p,\eps}}
^{p}|\mathscr I_{\eps}\big]\le \|h_{\eps}\|_{L_{p,\eps}}
^{p} $. A straightforward adaptation of the proof of Proposition \ref{teo:hB} allows to take $\overline A=+\infty$ if $X$ is a compound Poisson process. The constraint on $h_{\eps}\leq w_{\eps}$ for $1\leq p<2,$ is classical (see e.g. \cite{kerkyacharian1992density}). For instance, it is satisfied if $h_{\eps}$ is compactly supported.
Assumption  $\tfrac{v_{\Delta}(\eps)}{F_{\Delta}(\eps)}\leq \tfrac13$ is not restrictive, by means of Lemma \ref{lemma:momenti}, it holds for all $\eps\in(0,1]$
\begin{align}
\frac{v_\Delta(\eps)}{F_\Delta(\eps)}=\frac{v_\Delta(\eps)}{\Delta\lambda_\eps}+\frac{v_\Delta(\eps)}{\Delta}\Big(\frac{\Delta}{F_\Delta(\eps)}-\frac{1}{\lambda_\eps}\Big)\leq &\frac{v_\Delta(\eps)}{\Delta\lambda_\eps}+\frac{2\Delta\sigma^2(\eps)+2\Delta^2b_\nu(\eps)^2}{\Delta\eps^2}\Big(\frac{\Delta}{F_\Delta(\eps)}-\frac{1}{\lambda_\eps}\Big)\nonumber \\ &\xrightarrow[\Delta\to 0]{} 0.\label{eq:vF}
\end{align} Moreover, on the class 
$\mathscr L_{M,\alpha}$, $\alpha\in(0,1)$ and $M>0$ defined in \eqref{eq:Lalpha}, Theorems \ref{thm:vd1} and \ref{thm:Fd1}  permit to derive a bound in
$\frac{v_\Delta(\eps)}{F_\Delta(\eps)}\leq C_{\alpha,M,\eps} {\Delta },$  when $\eps$ is fixed.
A similar result can be obtained if $\alpha\in(1,2)$ assuming additionally that $f$ is symmetric: Theorems \ref{thm:vd2} and \ref{thm:Fd2bis} give a bound of the form $\frac{v_\Delta(\eps)}{F_\Delta(\eps)}\leq C_{\alpha,M,\eps} {\Delta ^{1/\alpha}},$  when $\eps$ is fixed. Finally, if $\alpha=1$, the same holds by replacing $\Delta$ with $\Delta\log(1/\Delta)$. 
\smallskip

To get unconditional bounds we introduce the following result.

 \begin{lemma}\label{lem:nB}
Let $F_\Delta(\varepsilon):=\PP(|X_\Delta|>\varepsilon)$. For all $r\geq 0$ we have 
\begin{align*}
\Big(\frac{3nF_{\Delta}(\eps)}{2}\Big)^{-r}\leq\E\big[\mathbf{n}(\varepsilon)^{-r}{\big|\nB\geq 1}\big]\leq 2\exp\big(\tfrac{-3nF_{\Delta}(\eps)}{32}\big)+\Big(\frac{nF_{\Delta}(\eps)}{2}\Big)^{-r}.
\end{align*}
\end{lemma}

Using Lemma \ref{lem:nB}, we remove the conditioning on $\mathscr I_{\eps}.$ The terms appearing in the following upper bound are discussed in the more general setting of Section \ref{sec:mainthmdiss} below.

\begin{cor}\label{cor:hB}
Fix $0<\eps\le 1$, assume that $f$ belongs to the class $ \mathscr F({s,p,q,\mathfrak{M}_\eps,A(\eps)})$ defined in \eqref{eq:Fbesov}, for some $1\leq q\leq \infty$, $ 1\leq p<\infty$ and $\eps<\overline{A}<\infty$. If $1\leq p<2$ suppose that $h_{\eps}(x)\leq w_{\eps}(x)$, $\forall x\in A(\eps)$ for some symmetric function $w_{\eps}\in L_{p/2}$.
 Let $r>s>\frac{1}{p}$ and let $\widehat h_{n,\varepsilon}$ be the wavelet estimator of 
 $h_{\eps} $ on $A(\eps)$, defined in \eqref{eq hat heps}.
For any $n\geq 1$ and $\Delta>0$ such that 
$n{F_{\Delta}(\eps)}\geq1$ and
$\tfrac{v_{\Delta}(\eps)}{F_{\Delta}(\eps)}\leq\tfrac13$, the following inequality holds.
For all $J\in \N$ and $p\geq 2$: 
\begin{align*}
\E\big[\|\hat h_{n,\eps}&(\{X_{i\Delta}-X_{(i-1)\Delta}\}_{i\in \mathscr I_{\eps}})-h_{\eps}\|_{L_{p,\eps}}
^{p}\big]\leq \|h_{\eps}\|_{L_{p,\eps}}^{p}\big(1-F_{\Delta}(\eps)\big)^{n}\\
&+ C\bigg\{{2^{2Jp}}\Big[\Big(\frac{v_{\Delta}(\eps){e^{-\lambda_\eps\Delta}}}{n F_{\Delta}(\eps)^{2}}\Big)^{\frac p2}+\Big(\frac{v_{\Delta}(\eps){e^{-\lambda_\eps\Delta}}}{F_{\Delta}(\eps)}\Big)^{p}\Big]\\
 &+2^{-Jsp}+2^{J\frac p2}\ell_{p,\eps}^{p/2}\big(n F_{\Delta}(\eps)\big)^{-p/2}+{\Big(\frac{2^J}{n F_\Delta(\eps)}\Big)^{(p-1)}}\mathbf{1}_{p\ge 2}+(e^{\lambda_\eps\Delta}\Delta\|f\|_{L_{p,\eps}})^{p}
 \\
 & \hspace{-1.5cm}+{2^{J(5\frac p2-1)}}\Big[(nF_{\Delta}(\eps))^{1-p} (\Delta \mu_p(\eps)+(\Delta\sigma^2(\eps))^{\frac p2})+(nF_{\Delta}(\eps))^{-\frac p2}\big(\sigma^{2}(\eps)\Delta\big)^{\frac p2}+(b_\nu(\varepsilon)\Delta)^p\Big]\bigg\},
 \end{align*}  
 for  $\ell_{p,\eps}:=\|h_{\eps}\|_{L_{p/2,\eps}}\mathbf{1}_{p\ge 2}+\|w_{\eps}\|_{L_{p/2,\eps}}\mathbf{1}_{1\le p<2}$ and $C>0$ depending only on $s$, $p$, $\|\Phi\|_{\infty}$, $\|\Phi'\|_{\infty}$, $\|\Phi\|_{p}$ and $\mathfrak M$.
\end{cor}

\section{Statistical properties of $\widehat f_{n,\varepsilon}$}\label{sec:mainthm}

\subsection{Main Theorem: general result \label{sec:mainthm1} }

Combining the results in Proposition \ref{teolambda} and Corollary \ref{cor:hB} we derive the following upper bound for the estimator $\widehat f_{n,\varepsilon}$ of the Lévy density $f$, when $\nu(\R) = \infty$. 

\begin{theorem}\label{thm:main}
Fix $0<\eps\le 1$, suppose $f$ belongs to the class $\mathscr F({s,p,q,\mathfrak{M}_\eps,A(\eps)})$ defined in \eqref{eq:Fbesov}, for some $1\leq q\leq \infty$, $ 1\leq p<\infty$ and $\eps<\overline{A}<\infty$. If $1\leq p<2$ suppose that $f_{\eps}(x)\leq \lambda_{\eps }w_{\eps}(x)$, $\forall x\in A(\eps)$ for some symmetric function $w_{\eps}\in L_{p/2}$.
 Let $r>s>\frac{1}{p}$ and let $\widehat f_{n,\varepsilon}$ be the estimator of 
 $f$ on $A(\eps)$, defined in \eqref{eqhatf}.
For any $n\geq 1$ and $\Delta>0$ such that $n{F_{\Delta}(\eps)}\geq1$ and
$\tfrac{v_{\Delta}(\eps)}{F_{\Delta}(\eps)}\leq\tfrac13$, the following inequality holds.
For all $J\in \N$, $p\geq 2$, there exists $C>0$ such that:
  \begin{align*}
  \big[\ell_{p,\varepsilon}\big(&\widehat f_{n,\varepsilon},f\big)\big]^p\leq C\bigg\{ \|f\|^{p}_{L_{p,\eps}}\big(1-F_{\Delta}(\eps)\big)^{n}+\Big(\frac{F_{\Delta}(\eps)}{n\Delta^2}\Big)^{\frac{p}{2}}+\Big|\lambda_\eps-\frac{F_{\Delta}(\eps)}{\Delta}\Big|^p
 \\ & +\lambda_{\eps}^{p}\Big\{{2^{2Jp}}\Big[\Big(\frac{v_{\Delta}(\eps){e^{-\lambda_\eps\Delta}}}{n F_{\Delta}(\eps)^{2}}\Big)^{p/2}+\Big(\frac{v_{\Delta}(\eps){e^{-\lambda_\eps\Delta}}}{F_{\Delta}(\eps)}\Big)^{p}\Big]\\
 &+2^{-Jsp}+2^{Jp/2}\ell_{p,\eps}^{p/2}\big(n F_{\Delta}(\eps)\big)^{-p/2}+{\bigg(\frac{2^J}{nF_\Delta(\eps)}\bigg)^{p-1}}\mathbf{1}_{p\ge 2}+(e^{\lambda_\eps\Delta}\Delta\|f\|_{L_{p,\eps}})^{p}
 \\
 & \hspace{-1cm}+{2^{J(5\frac p2-1)}}\Big[(nF_{\Delta}(\eps))^{1-p} \big(\Delta\mu_p(\eps)+(\Delta\sigma^2(\eps))^{\frac p2}\big)+(nF_{\Delta}(\eps))^{-\frac p2}\big(\sigma^{2}(\eps)\Delta\big)^{\frac p2}+(b_\nu(\varepsilon)\Delta)^p\Big]\Big\}\bigg\}
  \end{align*}
 where  $v_{\Delta}(\eps):=\PP(|M_{\Delta}(\varepsilon)+\Delta b_\nu(\varepsilon) |>\eps)$, $F_{\Delta}(\eps):=\PP(|X_{\Delta}|>\eps)$, $\sigma^2(\varepsilon):=\int_{|x|\leq \varepsilon}x^2\nu(dx)$, $\mu_p(\eps):=\int_{|x|\leq \eps}|x|^p\nu(dx)$,  $\ell_{p,\eps}:=\|h_{\eps}\|_{L_{p/2,\eps}}\mathbf{1}_{p\ge 2}+\|w_{\eps}\|_{L_{p/2,\eps}}\mathbf{1}_{1\le p<2}$ and $C$ depends on $s$, $p$, $\|\Phi\|_{\infty}$, $\|\Phi'\|_{\infty}$, $\|\Phi\|_{p}$ and $\mathfrak M$.
\end{theorem}

\subsubsection{Comments on the upper bound\label{sec:mainthmdiss}}
Theorem \ref{thm:main} gives an explicit upper bound for the $L_p$-risk restricted to the estimation set $A(\eps)$ for the estimation of the Lévy density $f=\frac{\nu(dx)}{dx}$ of a pure jump Lévy process as in \eqref{eq:defX}. Note that Lemma \ref{lemma:momenti} (see also \eqref{eq:vF}) ensures that $v_{\Delta}(\eps)/F_{\Delta}(\eps)$ tends to 0 when $\Delta\to0$. Therefore the upper bound of Theorem \ref{thm:main} tends to 0 whenever $\Delta\to 0$ and $n\Delta\to\infty$, under the assumption that the Lévy density $f$ has a Besov norm restricted to $A(\eps)$ that does not grow more than a constant times $\lambda_\eps$.

 We provide below a rough intuition of the different terms appearing in Theorem \ref{thm:main}. The estimation strategy relies on different approximations that entail four different sources of errors (points 3-4-5 hereafter are related to the estimation of $h_{\eps}$). 
\begin{enumerate}
\item
\textit{Controlling the presence of jumps}:  The term $\|f_{\eps}\|^{p}_{L_{p,\eps}}\big(1-F_{\Delta}(\eps)\big)^{n}$ provides a control of the risk when no jumps larger than $\eps$ are observed in a dataset. This term is bounded by $\|f_{\eps}\|^{p}_{L_{p,\eps}}e^{-nF_{\Delta}(\eps)}$ and tends to 0 exponentially fast under the assumption $nF_{\Delta}(\eps)\to\infty$ as $n\to\infty$. 
\item\textit{Estimation of $\lambda_\varepsilon$:} it leads to the error 
$\Big(\frac{F_{\Delta}(\eps)}{n\Delta^2}\Big)^{\frac{p}{2}}+\Big|\lambda_\eps-\frac{F_{\Delta}(\eps)}{\Delta}\Big|^p:=E_1.$
\item
\textit{Neglecting the event $\{|M_{\Delta}(\varepsilon)+\Delta b_\nu(\varepsilon)|>\eps\}$:} Considering that each time an increment $X_{\Delta}$ exceeds the threshold $\eps$ the associated Poisson process $N_{\Delta}(\eps)$ is nonzero leads to the error $$2^{2J}\bigg\{\sqrt{\frac{v_{\Delta}(\eps){e^{-\lambda_\eps\Delta}}}{n F_{\Delta}(\eps)^{2}}}+\frac{v_{\Delta}(\eps){e^{-\lambda_\eps\Delta}}}{F_{\Delta}(\eps)}\bigg\}\asymp 2^{2J}\frac{v_{\Delta}(\eps)}{F_{\Delta}(\eps)}:=E_{2}.$$ This error is unavoidable as we do not observe $M(\eps)$ and $Z(\eps)$ separately.
\item
\textit{Neglecting the presence of $M_{\Delta}(\varepsilon)+\Delta b_\nu(\varepsilon)$:} In \eqref{eq:M0} we ignore the convolution structure of the observations. This produces an error $E_{3}$ in $${\small 2^{J(\frac52-\frac1p)}\Big\{(nF_{\Delta}(\eps))^{-1+\frac1p} \big(\Delta \mu_p(\eps) +(\Delta\sigma^2(\eps))^{\frac p2}\big)+(nF_{\Delta}(\eps))^{-\frac12}\big(\sigma^{2}(\eps)\Delta\big)^{\frac12}+(b_\nu(\varepsilon)\Delta)^p\Big\}.}$$ It seems reasonable to neglect $M_{\Delta}(\varepsilon)+b_\nu(\varepsilon)\Delta$: the distribution of $M_{\Delta}(\varepsilon)$ is unknown. Even if we did know it, deconvolution methods essentially rely on spectral approaches which we meant to avoid.
\item
\textit{Estimation of the compound Poisson $Z(\varepsilon)$:} This estimation problem is solved in two steps. First, we neglect the event $\{N_{\Delta}(\varepsilon)\geq 2\}$ which generates the error:
$E_{4}:={e^{\lambda_\eps\Delta}}\Delta\|f\|_{L_{p,\eps}}.$ This could have been improved considering a corrected estimator as in \cite{Duval}, but it would have added even more heaviness in the final result. Second, for $2^J\lesssim n F_\Delta(\eps)$, we recover an estimation error that is classical for the density estimation pro\-blem from i.i.d. observations in $E_{5}:=2^{-Js}+2^{J/2}\ell_{p,\eps}^{1/2}\big(n F_{\Delta}(\eps)\big)^{-1/2}\hspace{-0.2cm}.$ 
\end{enumerate}

The rate of convergence is not explicit in terms of $\Delta$, it depends on the quantities $F_\Delta(\eps)$ and $v_\Delta(\eps)$ and therefore on the Lévy measure $\nu$. Consequently, we cannot say in general which one of the above error terms $E_{1}$, $E_{2}$, $E_{3}$, $E_{4}$ or $E_{5}$ is predominant. 

Moreover, the rate of convergence depends on the choice made for $J$. There is a bias term in $2^{-Jsp}$ (see $E_{5}$), where $s$ is the regularity of $f$, that decreases with $J$, whereas all the other ---variance type--- terms ($E_{2}$, $E_{3}$ and the second term of $E_{5}$) increase with $J$. Ideally, $J$ should be selected of the order of $J^{*}_{n}$ a minimizer of the upper bound of Theorem \ref{thm:main}. In the results below, from the idea that $E_{2}$ and $E_{3}$ are --under suitable assumptions-- remainder terms and should not intervene in the final rate, we  focus only on $E_{5}$ and select $J$ as a minimizer  of $E_{5}$. An adaptive procedure to select $J$ is discussed in Section \ref{sec:discussmain2}.

In Theorems \ref{thm:mainF} and \ref{thm:mainF2}   below, the rates of convergence for the $L_p$-risk $\ell_{p,\varepsilon}\big(\widehat f_{n,\varepsilon},f\big)$ are given under additional assumptions on $F_\Delta(\eps)$ and $v_\Delta(\eps)$, which are satisfied on the class $\mathscr{L}_{\alpha,M}$ defined in \eqref{eq:Lalpha}. Finally, Theorem \ref{thm:mainF3} deals with the case $\eps\to0$, it shows that our estimator is robust in this context.

\subsubsection{Example: Compound Poisson process\label{sec:mainthmex}} If $X$ is a compound Poisson process, we fix $\eps=0$. Then $\lambda:=\lambda_{0}<\infty$,  $F_{\Delta}(0)=1-e^{-\lambda_{0}\Delta}$  and $v_{\Delta}(0)=\mu_{p}(0)=\sigma^{2}(0)=b_{\nu}(0)=0$.  The bound given in Theorem \ref{thm:main} simplifies and choosing $J$ such that $2^{J}=(n\Delta)^{\frac{1}{2s+1}}$ we get
$  \big[\ell_{p,0}\big(\widehat f_{n,0},f\big)\big]^p\leq C\big\{
(n\Delta)^{-\frac{sp}{2s+1}}+\Delta^{p}\big\},$
where the first term is the optimal rate of convergence to estimate $p_{\Delta,0}$ from the observations $\bm D_{n,0}$ and the second term is a deterministic error due to the omission of the event that more than one jump may occur in an interval of length $\Delta$ (see also \cite{Duval}).

\subsection{An explicit rate under additional assumptions\label{sec:thmmain2}  }

\subsubsection{Regimes where $n\Delta^{2}\le 1$}

Without any specific assumption on the Lévy density, the rate of Theorem \ref{thm:main} is not explicit. This rate can be simplified under the additional assumptions \eqref{eq:assH1}, $n\Delta^{2}\le 1$ and $nv_{\Delta}(\eps)\le 1$  as well as choosing $J$ such that $2^{J}=(n\Delta)^{1/(2s+1)}$. The dimension $J$ is selected such that the estimator of $h_{\eps}$ is rate optimal from the observations $\mathbf{D}_{n,\eps}$. The following theorem is derived from Theorem \ref{thm:main} (see Appendix \ref{sec:app} for the proof).

\begin{theorem}\label{thm:mainF}
Let $X$ be a Lévy process as in \eqref{eq:defX} and let $\nu$ be a Lévy measure admitting a density $f$ with respect to the Lebesgue measure. Fix $0<\eps\le 1$, suppose that $f$ belongs to the class $\mathscr F_{\mathscr	H_{1}}:=\{f\in \mathscr	F({s,p,q,\mathfrak{M}_\eps,A(\eps)}),\ \lambda_\eps\ge 1, f\in \mathscr H_{1}\},$ 
 where $\mathscr F({s,p,q,\mathfrak{M}_\eps,A(\eps)})$ is defined in \eqref{eq:Fbesov}, for some $1\leq q\leq \infty$, $ 2\leq p<\infty$, and $\eps<\overline{A}<\infty$ and where  $\mathscr H_{1}:=\{f:\ \eqref{eq:assH1}\mbox{ is satisfied}\}$, for some $\delta>0$ and $c>0$.
 Let $r>s>\frac{1}{p}$ and let $\widehat f_{n,\varepsilon}$ be the estimator of 
 $f$ on $A(\eps)$, defined in \eqref{eqhatf}, let $n$ and $\Delta$ be such that $n\Delta\ge 6$, $n\Delta^{2}\le 1$, $\Delta\le \frac{1}{2c}\wedge\delta$ and $nv_{\Delta}(\eps)\le 1$. 
Then, for all $p\geq 2$, it holds that $\big[\ell_{p,\varepsilon}\big(\widehat f_{n,\varepsilon},f\big)\big]^p$ is bounded uniformly over $\mathscr F_{\mathscr	H_{1}}$ by
$$ C
\begin{cases}
(n\Delta)^{-\frac{sp}{2s+1}} & \mbox{if } s\ge\frac{3}{2}-\frac1p,
\\
\max\big\{(n\Delta)^{-\frac{sp}{2s+1}},\Delta^{1\wedge \frac p2}(n\Delta)^{-\frac{4s(p-1)-3p}{2(2s+1)}},(n\Delta)^{\frac{p(1-2s)}{(2s+1)}}\big\} & \mbox{else},
\end{cases}
$$ where $C$ is a constant depending on $\eps $, $s$, $p$, $\|\Phi\|_{\infty}$, $\|\Phi'\|_{\infty}$, $\|\Phi\|_{p}$ and $\mathfrak M$.
 \end{theorem}
 
If $1\le p<2$, the result of Theorem \ref{thm:mainF} follows from the case $p\ge 2$ and the H\"older inequality, as the set $A(\eps)$ is bounded and $f\in L_{2,\eps}\cap L_{\frac{2p}{p-2},\eps}$.

The quantity $v_\Delta(\eps)$ appearing in the hypotheses of Theorem \ref{thm:mainF} can be easily controlled on the class of subordinators satisfying \eqref{eq:assH1}, $\Delta\le \delta$, for instance for a Gamma process. Indeed, for any subordinator $X$ it holds 
$\PP(X_{\Delta}>\eps)=v_{\Delta}(\eps)e^{-\lambda_{\eps}\Delta}+1-e^{-\lambda_{\eps}\Delta},$  under Assumption \eqref{eq:assH1} it becomes $v_{\Delta}(\eps)=e^{\lambda_{\eps}\Delta}(F_{\Delta}(\eps)+e^{-\lambda_{\eps}\Delta}-1)=O\big(\Delta^{2}\big)$ as $\Delta\to 0$.
More generally, it is possible to give explicit upper bounds for $v_\Delta(\eps) $ on a larger class of Lévy processes than the subordinators, namely the class $\mathscr{L}_{M,\alpha}$, $\alpha\in(0,1)$ and $M>0$, on which it is possible to show that 
\begin{equation}\label{eq:vdiscu}
v_{\Delta}(\eps)\le C_{\alpha,M,\eps}\Delta^{2},
\end{equation} 
see Theorem \ref{thm:vd1}. The same result can be obtained on the class $\mathscr{L}_{M,\alpha}$, $\alpha\in(1,2)$, under the additional assumption of a symmetric Lévy density which is also $M\eps^{-(2+\alpha)}$-Lipschitz on the interval $(3/4\eps,5/4\eps)$, see Equation \eqref{eq:Mplus}.
Concerning the generality of Assumption \eqref{eq:assH1} we emphasize that it is always satisfied for Lévy densities in $\mathscr{L}_{M,\alpha}$, $\alpha\in(0,1)$ (see Theorem \ref{thm:Fd1}) and  for symmetric Lévy densities in $\mathscr{L}_{M,\alpha}$, $\alpha\in[1,2)$ that are $M\eps^{-(2+\alpha)}$-Lipschitz on the interval $(3/4\eps,5/4\eps)$, see Theorem \ref{thm:Fd2}.

  \subsubsection{An explicit rate for all regimes}

\textbf{Assumption $H_2(\beta,\delta,c)$}: $X$ is a Lévy process as in \eqref{eq:defX} such that
\begin{align}
\PP(|M_t(\eps)+tb_\nu(\eps)|>\eps)\le ct^\beta,\quad \forall 0<t\leq \delta,\ \forall 0<\eps\le 1,\label{eq:assH2}\tag{$H_{2}(\beta,\delta,c)$}
\end{align}
for some $\beta>1$ and where $c$ is a positive constant that does not depend on $t$. \smallskip

Note that \eqref{eq:assH2} is satisfied for $\beta=2$ if $f$ belongs to the class $\mathscr{L}_{\alpha,M}$, $\alpha\in(0,1)$ (see \eqref{eq:Lalpha} and Theorem \ref{thm:vd1}) or if $X$ is a subordinator such that \eqref{eq:assH1} is satisfied. If $X$ is a symmetric infinite variation process,
\eqref{eq:assH2} is satisfied for $\beta=1+1/\alpha$ if $f$ belongs to the class $\mathscr{L}_{\alpha,M}$, $\alpha\in[1,2)$ (see Theorem \ref{thm:vd2}). Observe that the exponent $1+1/\alpha$ cannot be improved in general, see Section 2.2 in \cite{DMnote}). However, if, additionally, $f$ is $M\eps^{-(2+\alpha)}$-Lipschitz on the interval $(3/4\eps,5/4\eps)$, the result can be improved with $\beta=2$ (see \eqref{eq:Mplus}).
The case $n\Delta^{2}\le 1$ being covered by Theorem \ref{thm:mainF}, we concentrate on the cases  $\beta\in(1,2)$ or $n\Delta^{2}>1.$

\begin{theorem}\label{thm:mainF2}
Let $X$ be a Lévy process as in \eqref{eq:defX} and let $\nu$ be a Lévy measure admitting a density $f$ with respect to the Lebesgue measure. Fix $0<\eps\le 1$, suppose that $f$ belongs to the class $\mathscr F_{\mathscr	H_{1},\mathscr	H_{2}}$ defined by  \begin{align*}
\mathscr F_{\mathscr	H_{1},\mathscr H_{2}}:=\{f\in \mathscr	F({s,p,q,\mathfrak{M}_\eps,A(\eps)}),\ \lambda_{\eps}\ge 1, f\in \mathscr H_{1}\cap\mathscr H_{2}\},\end{align*}
 with $\mathscr F({s,p,q,\mathfrak{M}_\eps,A(\eps)})$ defined in \eqref{eq:Fbesov}, for some $1\leq q\leq \infty$, $ 2\leq p<\infty$, and $\eps<\overline{A}<\infty$, $\mathscr H_{1}:=\{f,\ \eqref{eq:assH1}\mbox{ is satisfied}\}$ and $\mathscr H_{2}:=\{f,\ \eqref{eq:assH2}\mbox{ is satisfied}\}$ for some $c>0$, $\delta>0$ and $\beta>1$.  Let $r>s>\frac{1}{p}$ and let $\widehat f_{n,\varepsilon}$ be the estimator of 
 $f$ on $A(\eps)$, defined in \eqref{eqhatf}. Suppose that $n\ge 1$ and  that $\Delta$ is such that $\Delta\le \frac{1}{6c}\wedge\delta$.

Then, for all $p\geq 2$, the $L_p$-risk $ \big[\ell_{p,\varepsilon}\big(\widehat f_{n,\varepsilon},f\big)\big]^p$  is bounded uniformly over $\mathscr F_{\mathscr	H_{1},\mathscr	H_{2}}$ by
$$C\begin{cases}
\max\Big( (n\Delta)^{-\frac{sp}{2s+1}},v_{3}\Big)   &\text{if } n\Delta_{n}^2> 1\   \text{ and } \ \beta\ge 2\\
\max\Big((n\Delta)^{-\frac{sp}{2s+1}}, v_{2},v_{3},\tilde v_{3}\Big)   &\text{if } n\Delta_{n}^2> 1\   \text{ and } \ \beta\in(1, 2), \\
\max\Big((n\Delta)^{-\frac{sp}{2s+1}},v_{1},v_{2},\tilde v_{3}\Big)  &\text{if } n\Delta_{n}^2\le 1\   \text{ and } \ \beta\in(1,2),
  \end{cases}
$$  where $v_{1}:=(n\Delta)^{\frac{3p-2sp}{2(2s+1)}}\Delta^{\frac p2(\beta-1)}$, $v_{2}:=\Delta^{(\beta-1)p}(n\Delta)^{\frac{2p}{2s+1}},$ $v_{3}:=\Delta^{p}(n\Delta)^{\frac{5p-2}{2(2s+1)}}$ and $\tilde v_{3}:=(n\Delta)^{\frac{5p-2}{2(2s+1)}}(n\Delta)^{1-p}\Delta$, $C$ is a constant depending on $\eps $, $s$, $p$, {\small$\|\Phi\|_{\infty}$, $\|\Phi'\|_{\infty}$, $\|\Phi\|_{p}$ }and $\mathfrak M$. 

 \end{theorem}
 
From the proof of Theorem \ref{thm:mainF2}, we have the following relations: $v_{1}\le v_{2}$ if and only if $n\Delta^{\beta}\ge1$. Moreover, if $n\Delta^{2}\le 1$, then $(n\Delta)^{-\frac{sp}{2s+1}}\vee \tilde v_{3}=(n\Delta)^{-\frac{sp}{2s+1}}$ if $s\ge 3/2-1/p$ and $v_{1}\vee\tilde v_{3}=v_{1}$ if $s\ge1/2$. If $n\Delta^{2}>1$, $\beta\ge 2$ and $s\le 3/2-1/p$ then  $v^{*}\vee  v_{3}=v_{3}.$

 \subsubsection{Estimation in a neighborhood of the origin\label{sec:eps0}}
 
In the following result, we make explicit the dependency in $\eps$ of the bound for $f\in\mathscr	L_{M,\alpha},\ \alpha\in(0,2).$ This shows that $\hat f_{n,\eps}$ is robust when $\eps:=\eps_{n}\to0$ slowly with respect to $\Delta$.

\begin{theorem}\label{thm:mainF3}
Let $X$ be a Lévy process as in \eqref{eq:defX} and let $\nu$ be a Lévy measure admitting a density $f$ with respect to the Lebesgue measure. Fix $\eps\in(0,1]$ and suppose that $f\in\mathscr L_{M,\alpha}\cap \mathscr F({s,p,q,\mathfrak{M}_\eps,A(\eps)}) $ for some $\alpha\in(0,2)$ and for some $1\leq q\leq \infty$, $ 2\leq p<\infty$, and $\eps<\overline{A}<\infty$ (see \eqref{eq:Lalpha} and   \eqref{eq:Fbesov}). If $\alpha\in[1,2)$ suppose additionally that $f$ is symmetric and is $M\eps^{-(2+\alpha)}$-Lipschitz on the interval $(3/4\eps,5/4\eps)$.

 Let $r>s>{3}/{2}-{1}/{p}$ and let $\widehat f_{n,\varepsilon}$ be the estimator of 
 $f$ on $A(\eps)$, defined in \eqref{eqhatf}. Suppose that $n$, $\Delta$ and $\eps$ are such that $\lambda_{\eps}\ge 1,$ $\lambda_{\eps}\eps^{\alpha}\le 1$, $\eps^{1+\alpha-\frac2p}n\Delta\lambda_{\eps}^{2}\ge 1$, $n(\Delta\lambda_{\eps}\eps^{-2\alpha})^{2}\le 1$ and $\Delta\le \overline \Delta$, where $\overline \Delta$ is defined in \eqref{eq:Delta34}.

Then, for all $p\geq 2$, it holds 
$$  \big[\ell_{p,\varepsilon}\big(\widehat f_{n,\varepsilon},f\big)\big]^p\leq C\lambda_{\eps}^{p}(\eps^{1+\alpha-\frac2p}n\Delta\lambda_{\eps}^{2})^{-\frac{sp}{2s+1}},$$
 where $C$ is a constant depending on $\alpha$, $M$, $s$, $p$, $\|f\|_{L_{p,1}}$, $\|f\|_{L_{p/2,1}},\ \lambda_{1},\ \mu_{r}(1),\ r\in\{1,2,p\}$, $\|\Phi\|_{\infty}$, $\|\Phi'\|_{\infty}$, $\|\Phi\|_{p}$ and $\mathfrak M$.

 \end{theorem}

Note that for $p=2$ and $X$ an $\alpha$-stable processes, $\alpha\in(0,2)$, using $\lambda_{\eps}= C_{\alpha,\lambda_{1}} \eps^{-\alpha}$, we recover for $n,\Delta$ and $\eps$ selected as in Theorem \ref{thm:mainF3}  that $  \ell_{2,\varepsilon}\big(\widehat f_{n,\varepsilon},f\big)\leq C\lambda_{{\eps}}(n\Delta\lambda_{\eps})^{-\frac{s}{2s+1}}.$

\subsection{Discussion\label{sec:discussmain2}} 

\subsubsection{General comments}\label{subsec:comm}

Theorem \ref{thm:mainF} ensures that for any finite variation process of the form \eqref{eq:defX} --whose Lévy density is automatically  in $\mathscr	L_{1,M}$-- in regimes such that $n\Delta^{2}\le1$, our estimator attains the rate $(n\Delta)^{-sp/(2s+1)}$ uniformly over a Besov class of regularity $s\ge \frac32-\frac1p$. This rate is also attained for any symmetric infinite variation process such that its Lévy density is in $\mathscr L_{M,\alpha}$ for some $\alpha\in[1,2)$ and is $M\eps^{-2+\alpha}$-Lipschitz on the interval $(2/3\eps,5/4\eps)$. As mentioned earlier, the Lipschitz condition is satisfied if $f(x)=L(x)/x^{1+\alpha}$ for all $x\in[-2,2]\setminus \{0\}$ where $L$ is a bounded differentiable function with bounded derivative. This regularity condition on $f$ on a neighborhood of $\eps$  is different from requiring that $ f\in\mathscr F({s,p,q,\mathfrak{M}_\eps,A(\eps)}).$

This result generalizes \cite{MR2565560}: there, under the assumptions $p=2$, $n\Delta^{2}\le 1$ and $x^{2}f$ is in a Sobolev class of regularity $s>\frac12$, the same rate is attained for the estimation of $x^{2}f$ for a $L_{2}$ loss function. However, in dimension $d=1$, the hypothesis $x^{2}f$ is in a Sobolev space with regularity $s>\frac12$ implies that $\|x^{2}f\|_{\infty}<\infty$ and that $f\in \mathscr L_{1,M}$ for $M:=\|x^{2}f\|_{\infty}$ (see Theorem 4.2 in \cite{MR2424078}).  Therefore, we generalize \cite{MR2565560} to $L_{p}$ loss functions and to  symmetric infinite variation Lévy processes in $\mathscr L_{M,\alpha}$, $\alpha\in(1,2)$.

Moreover, the results of Theorem \ref{thm:mainF2} generalizes the latter  to $L_{p}, \ p\ge 1$ loss functions and without assumptions on $\Delta$ other that $\Delta\to0$. In these cases the rate of our procedure may  be slower than $(n\Delta)^{-sp/(2s+1)}$. 

Theorem \ref{thm:mainF3} shows that our estimator is robust when $\eps$ gets close to the critical value 0. This form of results was not studied in the literature. Unsurprisingly, the rate deteriorates as $\eps$ gets close to the critical value 0 as it gets multiplied with the increasing quantity $\lambda_{\eps}^{\frac{p}{2s+1}}\eps^{-\frac{sp}{2s+1}\big(1+\alpha-\frac2p\big)} $.

\subsubsection{Discussion on optimality}

The question whether  the upper bound of Theorem \ref{thm:main} is optimal remains open. When estimating $f\mathbf{1}_{A(\eps)}$, for $\eps$ fixed,
 we in fact estimate a Poisson measure. 
It is well known that the minimax rate of convergence over Besov classes with regularity $s$ for the estimation of the Lévy density of a compound Poisson process from the observation of $n$ of its increments sampled at rate $\Delta$ is $(n\Delta)^{-s/(2s+1)}$, see e.g. \cite{Duval}. 

However, in the present context, the data available to estimate the Poisson measure are not the increments of a compound Poisson process but those of a general Lévy process. This is a more challenging problem and in particular we can deduce that the rate $(n\Delta)^{-s/(2s+1)}$ cannot be improved on the class $\mathscr F(s,p,q,\mathfrak M_{\eps},A(\eps))$ (see also Section 4 of \cite{Lopez09}). This rate is obtained in Theorem \ref{thm:mainF} under the additional assumptions $n\Delta^{2}\le 1$, $s\ge \frac32-\frac1p$, \eqref{eq:assH1} with $\Delta\leq \delta$ and $nv_{\Delta}(\eps)\le 1$. Remark that for Lévy densities belonging to the nonparametric class $\mathscr{L}_{\alpha,M}$ if $\alpha\in (0,1)$ or for smooth symmetric Lévy densities in $\mathscr{L}_{\alpha,M}$ with $\alpha\in[1,2)$, the last two conditions are automatically satisfied as soon as $n\Delta^2\leq C_{M,\alpha,\eps}^{-1}$, where $C_{M,\alpha,\eps}$ is as in \eqref{eq:vdiscu}. Theorem \ref{thm:mainF3} shows that this rate is robust to small values of $\eps$.

Finally, we observe that there exist examples of Lévy processes of infinite variation and with non-smooth Lévy densities that do not satisfy Assumption \eqref{eq:assH1}, see e.g. \cite{marchal2009small}. Upper bounds for the quantities $F_\Delta(\eps)$ and $v_\Delta(\eps)$ are still known in this situation, see \cite{DMnote}. More precisely, for symmetric Lévy densities in $\mathscr{L}_{\alpha,M}$ with $\alpha\in(1,2)$, it holds
\begin{equation}\label{fveps}
|F_\Delta(\eps)-\lambda_\eps\Delta|\leq K_{M,\alpha,\eps}\Delta^{1+1/\alpha} \quad\text{and}\quad v_\Delta(\eps)\leq J_{M,\alpha,\eps} \Delta^{1+1/\alpha},
\end{equation}
for some constants $K_{M,\alpha,\eps}$ and $J_{M,\alpha,\eps}$ whose dependency in $M,\alpha$ and $\eps$ can be made explicit.
Using \eqref{fveps} in Theorem \ref{thm:main} leads to a rate that depends on $\alpha$ and which is slower than $(n\Delta)^{-s/(2s+1)}$. 
Nevertheless it is hard to say whether such a rate of convergence is  optimal.

\subsubsection{Adaptive selection procedure for $J$}

Motivated by the belief that the dominating term in the upper bound of Theorem \ref{thm:main} is $2^{-Js}+2^{J/2}\ell_{p,\eps}^{1/2}(n F_{\Delta}(\eps))^{-1/2}$, in Theorems \ref{thm:mainF} and \ref{thm:mainF2}  ($\eps$ is fixed) we selected $J_{n}^{*}$ such that $2^{J_{n}^{*}}=(n\Delta)^{1/(2s+1)}.$ A similar quantity depending on $\eps$ is considered in Theorem \ref{thm:mainF3}. However, such $J^{*}_{n}$ depends on the unknown regularity $s$ and is not feasible in practice. We propose here a data-driven procedure to select a suitable dimension $J$.

Usually, for wavelet type estimators, adaptation is achieved by thresholding techniques. However, the study of our estimator relies on the property that $(X_{i\Delta})_{i}\mapsto \hat\alpha_{J,k}\big((X_{i\Delta})_{i}\big)$ is differentiable, which is no longer true if we threshold the coefficients $\big(\hat\alpha_{J,k}\big)_{k}$. 

Ideally, $J_{n}$ should be selected of the order of $J^{*}_{n}$ a minimizer of the upper bound of Theorem \ref{thm:main} or equivalently a term realizing the trade-off between this bias term $2^{-Jsp}$ and all the variance terms depending on $J$ (see $E_{2},\ E_{3}$ and the second term of $E_{5}$ in Section \ref{subsec:comm}).
The following adaptive procedure to select $J_{n}$ from the observations is inspired from the Goldenshluger and Lepskii's method see e.g. \cite{goldenshluger2008universal}. Select $\hat J_{n}$ such that an estimator of the bias is of the order of the variance.
Denote by $\hat f_{n,\eps}:=\hat f_{n}(J)$ (see \eqref{eqhatf} and \eqref{eq hat heps}), and consider the data driven choice  \begin{align}\label{eq:Jhat}
\hat J_{n}:=\min\Big\{J\in\mathcal{J},\ \|\hat f_{n}(J)-\hat f_{n}(J')\|_{p}^{p}\leq \kappa V_{p}(J),\, J\leq J',\, J'\in \mathcal{J}\Big\},
\end{align} where $\mathcal{J}=\{1,\ldots, \lfloor \log(T)/\log(2)\rfloor\}$, $\kappa>0$ is a constant to be calibrated and $V_{p}(J)$ is an upper bound of the variance term appearing in the bound of Theorem \ref{thm:main}.  We do not seek for dimensions larger than $\lfloor \log(T)/\log(2)\rfloor$ as then the variance term no longer tends to 0 (see $E_{5}$ in Section \ref{subsec:comm}). However, the problem of the definition \eqref{eq:Jhat} is that it requires an explicit, sharp, upper bound for $V_{p}(J)$ which is not available without additional assumptions on the Lévy density $f$, $n$ and $\Delta$.  Under the Assumptions of Theorem \ref{thm:mainF}, we can set,
  \[V_{p}(J)=\frac{2^{Jp/2}}{(n\Delta)^{p/2}}+\frac{2^{2Jp}}{(n\Delta)^{p}}+\frac{2^{J(\frac{5p}{2}-1)}}{(n\Delta)^{p-1}}\Delta,\quad p\ge2.\]
 We do not investigate the question whether the resulting estimator $\hat f_{n}(\hat J_{n})$ satisfies the same upper bound --up to a numerical constant-- as $\hat f_{n}(J^{*}_{n})$. A key element in establishing such a bound would be the control of the deviations of $\hat J_{n}$ from $J^{*}_{n}$.

\subsection{Extensions\label{sec:brown}\label{sec:gen}}
If the Lévy process has a Brownian component, the estimator presented here applies and the results established can be generalized without technical difficulties. Let $\widetilde X$ be a Lévy process of the form:
\begin{align*}
 \widetilde X_t&= tb_\nu(\varepsilon)+\sigma W_{t}+M_t(\varepsilon)+Z_t(\varepsilon),
\end{align*} where $W$ is a standard Wiener process, independent of $M(\eps)$ and $Z(\eps)$.
This corresponds to consider $\tilde X_t=\sigma W_t+ X_t$ with $X$ as in \eqref{eq:defX}.

Similarly, consider the increments $\big(\widetilde X_{i\Delta}-\widetilde X_{(i-1)\Delta}, i\in\widetilde{\mathscr{I}}_{\eps}\big)$ where $\widetilde{ \mathscr{I}}_{\eps}:=\big\{i=1,\ldots,n:|\widetilde X_{i\Delta}-\widetilde X_{(i-1)\Delta}|>\eps\big\}$. Equations \eqref{eq:lambdahat}  and \eqref{eq hat heps} applied to these increments give estimators of $\lambda_{\eps}$ and $h_{\eps}$. Since $\sigma W_{\Delta}\overset{d}{=}\sigma\sqrt{\Delta}\mathcal{N}(0,1)$, in the asymptotic $\Delta\to 0$ approximation \eqref{eq:M0} still makes sense when applied to the increments of $\widetilde X$. 

Theorem \ref{teolambda} holds true in this setting, its proof remains unchanged by the additional Brownian part. However, the quantity $|\lambda_{\eps}-\frac{1}{\Delta}\PP(|\widetilde X_{\Delta}|>\eps)|$ needs to be handled differently in examples. Moreover, at the expense of small modifications in its proof, Corollary \ref{cor:hB} still holds after replacing in its statement
 $v_{\Delta}(\eps)$ by $v_{\Delta,\sigma}(\eps):=\PP\big(|\Delta b_{\nu}(\eps)+M_{\Delta}(\eps)+\sigma W_{\Delta}|>\eps\big)$ (this quantity agrees with the previous definition when $\sigma=0$). 

Combining those results,  Theorem \ref{thm:main} holds.
Applying Lemma \ref{lemma:momenti}, we recover that $\PP(|\tilde X_{\Delta}|>\eps)=O(\Delta)$ and $v_{\Delta,\sigma}(\eps)=o(\Delta)$ as $\Delta\to 0$, from which we obtain the consistency of the procedure. Moreover, Theorems \ref{thm:mainF}, \ref{thm:mainF2} and \ref{thm:mainF3} can also be obtained adapting \eqref{eq:assH2} accordingly, replacing $M_{t}(\eps)+tb_{\nu}(\eps)$ with $M_{t}(\eps)+tb_{\nu}(\eps)+\sigma W_{t}.$ Finally, as discussed in Section 2.5 of \cite{DMnote} the results of Appendix \ref{app:note} also hold in presence of a Gaussian component and Assumptions \eqref{eq:assH1} and \eqref{eq:assH2} are satisfied on the class $\mathscr L_{M,\alpha}$, $M>0,\ \alpha\in(0,2)$.\\

Another interesting extension would be to apply this methodology to estimate the Lévy density of an Itô semimartingale from high frequency observations as studied in \cite{hoffmann2017weak}. Then, for $\Delta$ small enough and under suitable assumptions, the contribution of the drift and the diffusive part of the process can be neglected and the above strategy should generalize. However, it would induce more technicalities in the proofs as the additional drift and diffusive part are not necessarily constant nor deterministic.

\section{Proofs\label{sec:prf}}
In the sequel, $C$ is a constant whose value may vary from line to line. Its dependencies may be given in indices.  Proofs of auxiliary lemmas are postponed to Appendix \ref{sec:app}.

\subsection{Proof of Theorem \ref{teolambda}}

Let $F_{\Delta}(\varepsilon):=\PP(|X_{\Delta}|> \eps)$ and $\widehat F_{\Delta}(\eps):=\frac{1}{n}\sum_{i=1}^{ n}\1_{(\varepsilon,\infty)} (|X_{i\Delta}-X_{(i-1)\Delta}|)$. The following holds
\begin{align}\label{eq:lambda}
 \E_{P_{n}}\Big[\big|\lambda_\eps-\widehat \lambda_{n,\eps}\big|^p\Big]&\leq2^p\bigg\{\Big|\lambda_\eps-\frac{F_{\Delta}(\eps)}{\Delta}\Big|^p+\frac{1}{\Delta^{p}}\E_{P_{n}}\Big[
 \big|F_{\Delta}(\eps)-\widehat F_{\Delta}(\eps)\big|^p\Big]\bigg\}.
 \end{align}
To control the second term in \eqref{eq:lambda}, we introduce the i.i.d. centered random variables $$U_i:=\frac{\1_{(\varepsilon,\infty)} (|X_{i\Delta}- X_{(i-1)\Delta}|)-F_{\Delta}(\eps)}{n},\quad i=1,\dots,n.$$  For $p\geq 2$, an application of the Rosenthal
inequality together with $\E\big[|U_i|^p\big]=O\big(\frac{F_\Delta(\varepsilon)}{n^p}\big)$ ensure the existence of a constant $C_p$ such that
$$\E_{P_{n}}\bigg[\Big|\sum_{i=1}^{ n}U_i\Big|^p\bigg]\leq C_p\bigg(n^{1-p}F_\Delta(\varepsilon)+\Big(\frac{F_{\Delta}(\eps)}{n}\Big)^{p/2}\bigg).$$ For $p\in[1,2)$, the Jensen inequality and the previous result for $p=2$ lead to
$$\E_{P_{n}}\bigg[\Big|\sum_{i=1}^{n}U_i\Big|^p\bigg]\leq \bigg(\E_{P_{n}}\Big[\Big|\sum_{i=1}^{ n}U_i\Big|^2\Big]\bigg)^{p/2}\leq \bigg(\frac{F_{\Delta}(\eps)}{n}\bigg)^{p/2}.$$

In the asymptotic $n\to\infty,$ using \eqref{eq:lambda}, we are only left to show that, for $p\geq 2$,
\begin{equation}
 \frac{1}{\Delta^p}\E_{P_n}\bigg[\bigg|\sum_{i=1}^n U_i\bigg|^p\bigg]=O\bigg(\Big(\frac{F_{\Delta}(\eps)}{n\Delta^2}\Big)^{p/2}\bigg). \label{eq:desiredasymptotic}
\end{equation}
An application of the Bernstein inequality (using that $|U_i|\leq n^{-1}$ and the fact that the variance $\V[U_i]\leq \frac{F_\Delta(\varepsilon)}{n^2}$) allows us to deduce that
$$\PP\bigg(\bigg|\sum_{i=1}^n U_i\bigg|\geq t\bigg)\leq 2\exp\bigg(-\frac{t^2n}{2F_\Delta(\varepsilon)+\frac{2t}{3}}\bigg).$$
Therefore,
\begin{align*}
 \E_{P_n}\bigg[\bigg|\sum_{i=1}^n U_i\bigg|^p\bigg]=p\int_0^\infty t^{p-1}\PP\bigg(\bigg|\sum_{i=1}^n U_i\bigg|\geq t\bigg)dt\leq2p\int_0^\infty t^{p-1}\exp\bigg(-\frac{t^2n}{2F_\Delta(\varepsilon)+\frac{2t}{3}}\bigg)dt.
\end{align*}
Observe that, for $t\leq \frac{3}{2}F_\Delta(\varepsilon)$, the denominator $2F_\Delta(\varepsilon)+\frac{2t}{3}$ is smaller than $3F_\Delta(\varepsilon)$ while for $t\geq \frac{3}{2}F_\Delta(\varepsilon)$ we have $2F_\Delta(\varepsilon)+\frac{2t}{3}\leq 2t$. It follows, after a change of variables, that
{\small\begin{align}
 \int_0^\infty t^{p-1}\exp\bigg(-\frac{t^2n}{2F_\Delta(\varepsilon)+\frac{2t}{3}}\bigg)dt&\leq  \frac{1}{2}\Big(\frac{3F_\Delta(\varepsilon)}{n}\Big)^{p/2}\Gamma\Big(\frac{p}{2}\Big) +
 \Big(\frac{2}{n}\Big)^p\Gamma\Big(p, \frac{nF_\Delta(\varepsilon)}{4}\Big),\label{eq:Gammap}
\end{align}}
where, $\Gamma(s,x) = \int_x^\infty x^{s-1}e^{-x}dx$ is the incomplete Gamma function and $\Gamma(s) = \Gamma(s,0)$ is the usual Gamma function. To conclude, we use the classical estimate for the incomplete Gamma function for $|x| \to \infty$:
$
\Gamma(s,x) \approx x^{s-1} e^{-x} \big(1 + \frac{s-1}{x} + O({x^{-2}})\big).
$
When \eqref{eq:Gammap} is divided by $\Delta^p$, it is asymptotically $O\big({(n\Delta)^{-p}} e^{-nF_\Delta(\varepsilon)}\big)$, which goes to 0 faster than \eqref{eq:desiredasymptotic}.\hfill$\Box$

\subsection{Proof of Proposition \ref{teo:hB}\label{sec:prfth2}}

\paragraph*{Preliminary}  As the proof is lengthy, we enlighten here the main difficulties arising from the fact that the estimator $\hat h_{n,\eps}$ uses the observations $\bm D_{n,\varepsilon}$, i.e. $\widehat h_{n,\varepsilon}=\widehat h_{n,\varepsilon}(\bm D_{n,\varepsilon})$.
As $\mathscr{I}_\eps\ne \emptyset$ it holds $\nB\geq 1$. \begin{enumerate}\item The cardinality of $\bm D_{n,\varepsilon}$ is $\bm n(\varepsilon)$ that is random. That is why in Proposition \ref{teo:hB} we study the risk of this estimator conditionally on $\mathscr I_{\eps}.$ 

\item  An observation of $\bm D_{n,\eps}$ is not a realization of $h_{\eps}$. Indeed, an increment of $Z(\varepsilon)$ does not necessarily correspond to one jump, whose density is $h_\varepsilon$. More demanding, the presence of the small jumps $M(\eps)$ needs to be taken into account. To do so we split the sample $\bm{D}_{n,\eps}$ in two parts according to the presence or absence of jumps in the Poisson part. On the subsample where the Poisson part is nonzero, we make an expansion at order 1 and we neglect the presence of the small jumps. \end{enumerate}

\paragraph*{Expansion of $\widehat h_{n,\eps}$} Consider $\bm D_{n,\eps}=\{X_{i\Delta}-X_{(i-1)\Delta},\ i\in \mathscr{I}_{\eps}\}$ the increments larger than $\eps$. Recall that, for each $i$, we have 
 $$X_{i\Delta}-X_{(i-1)\Delta}=\Delta b_\nu(\varepsilon)+M_{i\Delta}(\eps)-M_{(i-1)\Delta}(\eps) +Z_{i\Delta}(\eps)-Z_{(i-1)\Delta}(\eps).$$
 We split the sample as follows:
$
 \mathscr K_{\eps}:=\{i\in\mathscr I_{\eps},\ Z_{i\Delta}(\eps)-Z_{(i-1)\Delta}(\eps)\ne0\}$ and $ \mathscr K_{\eps}^{c}:= \mathscr I_{\eps}\setminus \mathscr K_{\eps}.$
Denote by $\tnB$ the cardinality of $\mathscr K_{\eps}$. To avoid cumbersomeness, in the remainder of the proof we write $M$ instead of $M(\eps)$ and $Z$ instead of $Z(\eps)$. 
 Recall that $\Phi_{Jk}(x)=2^{\frac J2}\Phi(2^{J}x-k)$. Using that $\Phi$ is continuously differentiable we can write, $\forall k\in \Lambda_{J}$,
\begin{align*}
\widehat	\alpha_{J,k}
&=\frac{1}{\nB}\sum_{i\in \mathscr K_{\eps}}\big\{\Phi_{Jk}(Z_{i\Delta}-Z_{(i-1)\Delta})+2^{3J/2}(M_{i\Delta}-M_{(i-1)\Delta}+b_\nu(\varepsilon)\Delta)\Phi'(2^{J}\eta_{i}-k)\big\}\\
&\quad +\frac{1}{\nB}\sum_{i\in \mathscr K_{\eps}^{c}}\Phi_{Jk}(X_{i\Delta}-X_{(i-1)\Delta}),
\end{align*} where $\eta_{i}\in [\min\{Z_{i\Delta}-Z_{(i-1)\Delta},X_{i\Delta}-X_{(i-1)\Delta}\},\max\{Z_{i\Delta}-Z_{(i-1)\Delta},X_{i\Delta}-X_{(i-1)\Delta}\}]$. It follows that
\begin{align*}
\hat h_{n,\eps}(x, \{X_{i\Delta}-&X_{(i-1)\Delta}\}_{i\in \mathscr I_{\eps}})=\sum_{k\in\Lambda_{J}}\widehat\alpha_{J,k}\Phi_{Jk}(x) :=\frac{\tnB}{\nB}\tilde h_{n,\eps}(x, \{Z_{i\Delta}-Z_{(i-1)\Delta}\}_{i\in \mathscr K_{\eps}})\\
&\ \ +\frac{2^{3J/2}}{\nB}\sum_{i\in \mathscr K_{\eps}}(M_{i\Delta}-M_{(i-1)\Delta} +b_\nu(\varepsilon)\Delta)\sum_{k\in\Lambda_{J}}\Phi'(2^{J}\eta_{i}-k)\Phi_{Jk}(x)\\
&\ \ +\frac{1}{\nB}\sum_{i\in\mathscr K_{\eps}^{c}}\sum_{k\in\Lambda_{J}}\Phi_{Jk}(M_{i\Delta}-M_{(i-1)\Delta}+b_\nu(\varepsilon)\Delta)\Phi_{Jk}(x),\end{align*} where conditional on $\mathscr K_{\eps}$, $\tilde h_{n,\eps}(\{Z_{i\Delta}-Z_{(i-1)\Delta}\}_{i\in \mathscr K_{\eps}})$ is the linear wavelet estimator of $\peps$ defined in \eqref{eq:peps} from $\tnB$ direct measurements. Explicitly, it is defined as follows
\begin{align}\label{eq:tildealphaJk}
\tilde h_{n,\eps}(x,\{Z_{i\Delta}-Z_{(i-1)\Delta}\}_{i\in \mathscr K_{\eps}})&=\sum_{k\in\Lambda_{J}}\widetilde\alpha_{J,k}\Phi_{Jk}(x),
\end{align} 
where $\widetilde\alpha_{J,k}=\frac{1}{\tnB}\sum_{i\in\mathscr K_{\eps}}\Phi_{Jk}(Z_{i\Delta}-Z_{(i-1)\Delta})$.
This is not an estimator as both $\mathscr K_{\eps}$ and $\{Z_{i\Delta}-Z_{(i-1)\Delta}\}_{i\in \mathscr K_{\eps}}$ are not observed. However, $\widetilde\alpha_{J,k}$ approximates the quantity \begin{align}
\label{eq:alphaJk}\alpha_{J,k}:=\int_{A(\eps)}\Phi_{Jk}(x)\peps(x)dx.
\end{align}

\paragraph*{Decomposition of the $L_{p,\eps}$ loss} Taking the $L_{p,\eps}$ norm and applying the triangle inequality we get
\begin{align}
\|\hat h_{n,\eps}&(\{X_{i\Delta}-X_{(i-1)\Delta}\}_{i\in \mathscr I_{\eps}})-h_{\eps}\|_{L_{p,\eps}}
^{p}\hspace{-0.1cm}\leq C_{p}\bigg\{ \Big(\frac{\tnB}{\nB}\Big)^{p}\|\tilde h_{n,\eps}(\{Z_{i\Delta}-Z_{(i-1)\Delta}\}_{i\in \mathscr K_{\eps}})-h_{\eps}\|_{L_{p,\eps}}^{p}\nonumber\\
 & \hspace{1cm}+\Big(1-\frac{\tnB}{\nB}\Big)^{p}\|h_{\eps}\|_{L_{p,\eps}}^{p}\nonumber\\
  & \hspace{1cm}+\frac{2^{3Jp/2}}{\nB^{p}}\int_{A(\eps)}\Big|\sum_{i\in\mathscr K_{\eps}}(M_{i\Delta}-M_{(i-1)\Delta}+b_\nu(\varepsilon)\Delta)\sum_{k\in\Lambda_{J}}\Phi'(2^{J}\eta_{i}-k)\Phi_{Jk}(x)\Big|^{p}dx\nonumber
  \\ &\hspace{1cm}+\frac{1}{\nB^{p}}\int_{A(\eps)}\Big|\sum_{i\in\mathscr K_{\eps}^{c}}\sum_{k\in\Lambda_{J}}\Phi_{Jk}(M_{i\Delta}-M_{(i-1)\Delta}+\Delta b_\nu(\varepsilon))\Phi_{Jk}(x)\Big|^{p}dx
\bigg\}\nonumber\\
&\hspace{1.1cm}= C_{p}\big\{T_{1}+T_{2}+T_{3}+T_{4}\big\}.\label{eq:T}
\end{align}
After taking expectation conditionally on $\mathscr I_{\eps}$ and $\mathscr K_{\eps}$, we bound each term separately. 

\begin{rem}\label{rem:T1}If $X$ is a compound Poisson process and we take $\eps=0$, then $\hat h_{n,\eps}=\tilde h_{n,\eps}$ (and $\bm{n}(0)=\tilde{\bm{n}}(0)$) and $T_{2}=T_{3}=T_{4}=0$.\end{rem}

\paragraph*{Control of $T_{1}$}  We have
\begin{align*}
\|\tilde h_{n,\eps}(\{Z_{i\Delta}-Z_{(i-1)\Delta}\}_{i\in \mathscr K_{\eps}})-h_{\eps}\|_{L_{p,\eps}}^{p}&\leq C_{p}\big\{ \|\tilde h_{n,\eps}(\{Z_{i\Delta}-Z_{(i-1)\Delta}\}_{i\in \mathscr K_{\eps}})-{\peps}\|_{L_{p,\eps}}^{p}\\ &\hspace{1cm}+\|\peps-h_{\eps}\|^{p}_{L_{p,\eps}}\big\}
=:C_{p}(T_{5}+T_{6}).
\end{align*} 
The deterministic term $T_{6}$ is bounded using Lemma \ref{lem:hpeps} by $(2\Delta e^{\lambda_\eps\Delta}\|f\|_{L_{p,\eps}})^{p}$.
Taking the expectation conditionally on $\mathscr I_{\eps}$ and $\mathscr K_{\eps}$ of $T_{5}$, we recover the linear wavelet estimator of $\peps$ studied by Kerkyacharian and Picard \cite{kerkyacharian1992density} (see their Theorem 2). For the sake of completeness we reproduce the main steps of their proof.
The control of the bias is the same as in  \cite{kerkyacharian1992density}. Noticing that Lemma \ref{lem:hepsBesov} implies $\peps \in\mathscr F\big(s,p,q, \mathfrak	M,A(\eps)\big)$ (see Lemma 5.1 in \cite{Duval}), we get
\begin{align*}
\E\big[T_{5}|\mathscr I_{\eps}, \mathscr K_{\eps}\big]&\leq C_{p}\bigg\{2^{-Jsp}{\mathfrak M}^{p}+2^{J(p/2-1)}\sum_{k\in\Lambda_{J}}\E[|\widetilde\alpha_{J,k}-\alpha_{J,k}|^{p}|\mathscr I_{\eps}, \mathscr K_{\eps}]\bigg\},
\end{align*} where $\tilde\alpha_{J,k}$ and $\alpha_{J,k}$ are defined in \eqref{eq:tildealphaJk} and \eqref{eq:alphaJk}. First consider the case $p\geq 2$. We start by observing that
\begin{align*}
 &\E\bigg[\bigg|\frac{1}{\tnB}\sum_{i\in \mathscr K_\eps}\Phi_{Jk}(Z_{i\Delta}-Z_{(i-1)\Delta})-\int_{A(\varepsilon)}\Phi_{Jk}(x)p_{\Delta,\eps}(x)\bigg|\big|\mathscr I_\eps, \mathscr K_\eps\bigg]=\\&\sum_{I\subset \{1,\dots,n\}}\1_{\{I=\mathscr K_\eps\}}\E\bigg[\bigg|\frac{1}{|I|}\sum_{i\in I}\Phi_{Jk}(Z_{i\Delta}-Z_{(i-1)\Delta})-\int_{A(\varepsilon)}\Phi_{Jk}(x)p_{\Delta,\eps}(x)\bigg|\bigg],
\end{align*}
where $|I|$ denotes the cardinality of the set $I$. To bound the last term we apply the inequality of  Bretagnolle and Huber
to the i.i.d. centered random variables $\big(\Phi_{Jk}(Z_{i\Delta}-Z_{(i-1)\Delta})-\E[\Phi_{Jk}(Z_{i\Delta}-Z_{(i-1)\Delta})]\big)_{i\in I}$ bounded by $2^{J/2+1}\|\Phi\|_{\infty}$, conditional to $\{\mathscr K_\eps=I\}$. We obtain the bound on the previous term
\begin{align*}
C_{p}
\sum_{k\in\Lambda_{J}}\bigg\{&\frac{1}{|I|^{p/2}}\Big[2^{J}\int_{A(\eps)}\Phi(2^{J}x-k)^{2}\peps(x)dx\Big]^{p/2}+\\& \ \mathbf{1}_{p\ge 2}\frac{2^{(p-2)(J/2+1)}\|\Phi\|_{\infty}^{p-2}}{|I|^{p-1}}\int_{A(\eps)}2^{J}\Phi(2^{J}x-k)^{2}\peps(x)dx\bigg\}.
 \end{align*}
Therefore, we get
\begin{align*}
\sum_{k\in\Lambda_{J}}\E[|\widetilde\alpha_{J,k}-\alpha_{J,k}|^{p}|\mathscr I_{\eps},\mathscr K_{\eps}]&\leq C_{p}
\sum_{k\in\Lambda_{J}}\bigg\{\frac{1}{\tnB^{p/2}}\Big[2^{J}\int_{A(\eps)}\Phi(2^{J}x-k)^{2}\peps(x)dx\Big]^{p/2}\\&+\mathbf{1}_{p\ge 2}\frac{2^{(p-2)(J/2+1)}\|\Phi\|_{\infty}^{p-2}}{\tnB^{p-1}}\int_{A(\eps)}2^{J}\Phi(2^{J}x-k)^{2}\peps(x)dx\bigg\},
\end{align*}
where, as developed in  \cite{kerkyacharian1992density},
\begin{align*}&\sum_{k\in\Lambda_{J}}\Big[\int_{A(\eps)}2^{J}\Phi(2^{J}x-k)^{2}\peps(x)dx\Big]^{p/2}\leq {\mathfrak M}2^{J}\|\peps\|_{L_{p/2,\eps}}^{p/2}\le{\mathfrak M}2^{J}\|h_{\eps}\|_{L_{p/2,\eps}}^{p/2}\\ \mbox{and}\quad &\sum_{k\in\Lambda_{J}} \int_{A(\eps)}2^{J}\Phi(2^{J}x-k)^{2}\peps(x)dx\leq {\mathfrak M}2^{J}.
\end{align*}
 We can then conclude that
 \begin{align*}
\sum_{k\in\Lambda_{J}}\E[|\widetilde\alpha_{J,k}-\alpha_{J,k}|^{p}|\mathscr I_{\eps},\mathscr K_{\eps}]&\leq C_{p}{\mathfrak M}
\bigg\{\frac{2^{J}\|h_{\eps}\|^{p/2}_{L_{p/2,\eps}}}{{\tnB^{p/2}}}+\mathbf{1}_{p\ge 2}\frac{2^{Jp/2}\|\Phi\|_{\infty}^{p-2}}{\tnB^{p-1}}\bigg\}.
\end{align*} 
Plugging this last inequality in $T_{5}$ we obtain
$$\E[T_{5}|\mathscr I_{\eps},\mathscr K_{\eps}]\leq C\bigg\{2^{-Jsp}{\mathfrak M}^{p}+{\mathfrak M}\Big(\Big(\frac{2^{J}\|h_{\eps}\|_{L_{p/2,\eps}}}{\tnB}\Big)^{p/2}+\mathbf{1}_{p\ge 2}\Big(\frac{2^{J}}{\tnB}\Big)^{p-1}\Big)\bigg\},$$ where $C$ is a constant depending on $p$ and $\|\Phi\|_{\infty}$. Gathering all terms we get, for $p\geq2$, 
\begin{align*}
\E\Big[&\|\hat h_{n,\eps}(\{Z_{i\Delta}-Z_{(i-1)\Delta}\}_{i\in\mathscr K_{\eps}})-h_{\eps}\|^{p}_{L_{p,\eps}}|\mathscr I_{\eps},\mathscr K_{\eps}\Big]\nonumber\\& \leq C\bigg\{\Big(2^{-Js}{\mathfrak M}\Big)^{p}+{\mathfrak M}\Big(\Big(\frac{2^{J}{\|h_{\eps}\|_{L_{p/2,\eps}}}}{\tnB}\Big)^{p/2}+\mathbf{1}_{p\ge 2}{\Big(\frac{2^{J}}{\tnB}\Big)^{p-1}\Big)}+(e^{\Delta\lambda_\eps}\Delta\|f\|_{L_{p,\eps}})^{p}\bigg\},
\end{align*}where $C$ depends on $p$ and $\ \|\Phi\|_{\infty}$.
For $p\in [1,2)$, together with the additional assumption of $h_{\eps}$, following the lines of the proof of Theorem 2 of Kerkyacharian and Picard \cite{kerkyacharian1992density} we obtain the same bound as above replacing ${\|h_{\eps}\|_{L_{p/2,\eps}}}$ with $\| w_{\eps}\|_{L_{p/2,\eps}}$.
Finally, using  $\tnB\leq \nB$we have established for $p\ge 1$ that 
{\small\begin{align}
\label{eq:T1}\E\Big[T_{1}|\mathscr I_{\eps},\mathscr K_{\eps}\Big] \leq C\Big\{\Big(2^{-Js}{\mathfrak M}\Big)^{p}+{\mathfrak M}\Big(\Big(\frac{2^{J}{\ell_{p,\eps}}}{\tnB}\Big)^{p/2}+\mathbf{1}_{p\ge 2}{\Big(\frac{2^{J}}{\tnB}\Big)^{p-1}\Big)}+(e^{\Delta\lambda_\eps}\Delta\|f\|_{L_{p,\eps}})^{p}\Big\},
\end{align}}where $C$ depends on $p$ and $\ \|\Phi\|_{\infty}$ and $\ell_{p,\eps}:={\|h_{\eps}\|_{L_{p/2,\eps}}\mathbf{1}_{p\ge 2}+\|w_{\eps}\|_{L_{p/2,\eps}\mathbf{1}_{1\le p<2}}}$. 

Note that taking $J$ such that $2^{J}=\tnB^{\frac{1}{2s+1}}$ we have, uniformly over $\mathscr F(s,p,q,{\mathfrak M},A(\eps))$, an upper bound in $\tnB^{-s/(2s+1)}$ for the estimation of $h_{\eps}$, which is the optimal rate of convergence for a density from $\tnB$ direct independent observations (see  \cite{kerkyacharian1992density}). 
Moreover, we did not use that $A(\eps)$ is bounded to control this quantity, it was possible to have $\overline{A}=\infty$. 

\paragraph*{Control of $T_{3}$}  Using the fact that $\Phi'$ is compactly supported, we get
{\small\begin{align*}
\E[T_{3}&|\mathscr I_{\eps},\mathscr K_{\eps}]\leq\frac{2^{\frac{3Jp}{2}} \|\Phi'\|^{p}_{\infty}}{\nB^{p}}\int_{\R}\Big|\sum_{k\in\Lambda_{J}}\Phi(x-k)\Big|^{p}\frac{dx}{2^{J}}\E\Big[\Big|\sum_{i\in\mathscr K_{\eps}}\big(M_{i\Delta}-M_{(i-1)\Delta}+b_\nu(\varepsilon)\Delta\big)\Big|^{p}|\mathscr I_{\eps},\mathscr K_{\eps}\Big].
\end{align*}}
Furthermore, we use the following upper bound for the last term in the expression above:
\begin{align*}
\E\bigg[\Big|\sum_{i\in\mathscr K_{\eps}}&\big(M_{i\Delta}-M_{(i-1)\Delta}+b_\nu(\varepsilon)\Delta\big)\Big|^{p}|\mathscr I_{\eps},\mathscr K_{\eps}\bigg]\\
&\leq C_p\bigg\{\E\bigg[\Big|\sum_{i\in\mathscr K_{\eps}}\big(M_{i\Delta}-M_{(i-1)\Delta}\big)\Big|^{p}|\mathscr I_{\eps},\mathscr K_{\eps}\bigg]+\big(\tnB b_\nu(\varepsilon)\Delta\big)^p
\bigg\}.
\end{align*}
From the Rosenthal inequality conditional on $\mathscr I_{\eps}$ and $\mathscr K_{\eps}$ we derive for $p\geq2$
\begin{align*}\E\Big[\Big|\sum_{i\in\mathscr K_{\eps}}\big(M_{i\Delta}-M_{(i-1)\Delta}\big)\Big|^{p}|\mathscr I_{\eps},\mathscr K_{\eps}\Big]
&\leq C_{p}\Big\{\tnB \E[|M_{\Delta}|^{p}]+\big(\tnB\E[M_{\Delta}^{2}]\big)^{\frac{p}{2}}\Big\}.\end{align*} 
Observe that $\E[M_\Delta^2]=\Delta\sigma^2(\eps)$. There exists a constant $C_p'$, only depending on $p$, such that 
$\E[|M_\Delta|^p]\leq \Delta\int_{|x|\leq \eps}|x|^p\nu(dx)+C_p'(\Delta\sigma^2(\eps))^{p/2}.$Set $\mu_p(\eps)=\int_{|x|\leq \eps}|x|^p\nu(dx)$.

For $p\in[1,2)$ we obtain the same result using the Jensen inequality and the latter inequality with $p=2$. Next, \begin{align*}
\int_{\R}\Big|\sum_{k\in\Lambda_{J}}\Phi(x-k)\Big|^{p}\frac{dx}{2^{J}}&\leq2^{-J}|\Lambda_{J}|^{p}\|\Phi\|_{p}^{p}.
\end{align*} 
As $\Phi$ is compactly supported, and since we estimate $h_{\eps}$ on a set bounded by $\overline{A}$, for every $j\geq 0$, the set $\Lambda_J$ has cardinality bounded by $|\Lambda_{J}|\leq  C 2^{J},$ where $C$ depends on the support of $\Phi$ and $\overline{A}$.
It follows that,  
\begin{align}
\E[T_{3}|\mathscr I_{\eps},\mathscr K_{\eps}]\leq C& \|\Phi'\|^{p}_{\infty}\|\Phi\|_{p}^{p}{2^{J(5p/2-1)}}\bigg\{\tnB\nB^{-p}(\Delta \mu_p(\eps)+(\Delta\sigma^2(\eps))^{p/2}) \nonumber\\ &+\nB^{-p}\big(\tnB\sigma^{2}(\eps)\Delta\big)^{p/2}+\Big(\frac{\tnB b_\nu(\varepsilon)\Delta}{\nB}\Big)^p\bigg\}.\label{eq:T3}\end{align}

\paragraph*{Control of $T_{4}$} Similarly, for the last term we have\begin{align}
\label{eq:T4}\E[T_{4}|\mathscr I_{\eps},\mathscr K_{\eps}]&\leq{2^{J(2p-1)}}\|\Phi\|^{p}_{\infty}\|\Phi\|_{p}^{p}\Big(1-\frac\tnB\nB\Big)^{p}.
\end{align}

\paragraph*{Deconditioning on $\mathscr K_{\eps}$} Substituting \eqref{eq:T1}, \eqref{eq:T3} and \eqref{eq:T4} into \eqref{eq:T}, and noticing that $T_{2}$ is negligible compared to $\E(T_{4}|\mathscr I_{\eps},\mathscr K_{\eps})$, we obtain 
\begin{align*}
\E\big[&\|\hat h_{n,\eps}(\{X_{i\Delta}-X_{(i-1)\Delta}\}_{i\in \mathscr I_{\eps}})-h_{\eps}\|_{L_{p,\eps}}
^{p}|\mathscr I_{\eps},\mathscr K_{\eps}\big]\\ &\leq C\bigg\{{2^{2Jp}}\Big(1-\frac{\tnB}{\nB}\Big)^{p}+2^{-Jsp}+\Big(\frac{2^{J}{\ell_{p,\eps}}}{\tnB}\Big)^{p/2}+\Big(\frac{2^{J}}{\tnB}\Big)^{p-1}\mathbf{1}_{p\ge2}+(e^{\Delta\lambda_\eps}\Delta\|f\|_{L_{p,\eps}})^{p}
 \\
 &\quad+{2^{J(3p/2-1)}}|\Lambda_{J}|^{p}\Big[\tnB\nB^{-p} (\Delta \mu_p(\eps)+(\Delta\sigma^2(\eps))^{p/2})\\ &\quad+\nB^{-p}\big(\tnB\sigma^{2}(\eps)\Delta\big)^{p/2}+\Big(\frac{\tnB b_\nu(\varepsilon)\Delta}{\nB}\Big)^p\Big]\bigg\},\end{align*}
where $C$ depends on $s$, $p$, $\|\Phi\|_{\infty}$, $\|\Phi'\|_{\infty}$, $\|\Phi\|_{p}$ and $\mathfrak M$. To remove the conditional expectation on $\mathscr K_{\eps}$ we apply the following lemma, whose proof is postponed in the appendix.

\begin{lemma}\label{lem:tnB} Let $v_{\Delta}(\eps)=\PP(|M_{\Delta}(\eps)+\Delta b_\nu(\varepsilon) |>\eps)$ and $F_{\Delta}(\eps)=\PP(|X_{\Delta}|> \eps)$. If $\frac{v_{\Delta}(\eps)}{F_{\Delta}(\eps)}\leq \frac13$ and $\nB\geq 1$ , then for all $r\geq 0$, there exists a constant $C$ depending on $r$ such that
\begin{align*}
\E\big[\tnB^{-r}\big|\mathscr I_{\eps}\big]&\leq C\nB^{-r}\\
\E\big[(\nB-\tnB)^{r}\big|\mathscr I_{\eps}\big]&\leq C\Big\{\Big(\nB\frac{v_{\Delta}(\eps)e^{-\lambda_{\eps}\Delta}}{ F_{\Delta}(\eps)}\Big)^{r/2}+\Big(\nB\frac{v_{\Delta}(\eps)e^{-\lambda_{\eps}\Delta}}{F_{\Delta}(\eps)}\Big)^{r}\Big\}.
\end{align*}   \end{lemma}
Finally, using Lemma \ref{lem:tnB} and that $\tnB\leq \nB$, 
we complete the proof. \hfill$\Box$

\subsection{Proof of Theorem \ref{thm:main}}
 Theorem \ref{thm:main}  is a consequence of Theorem \ref{teolambda} and Corollary \ref{cor:hB}.  For all $0<\eps\le 1$, we decompose $\ell_{p,\varepsilon}(\widehat f_{n,\varepsilon},f)$ as follows:
 \begin{align*}
  \big[\ell_{p,\varepsilon}\big(\widehat f_{n,\varepsilon},f\big)\big]^p
&\leq 2^{p-1}\E_{P_n}\Big[\big|\widehat\lambda_{n,\varepsilon}-\lambda_\varepsilon\big|^p\Big]\|h_\varepsilon\|_{L_{p,\varepsilon}}^p+2^{p-1}\E_{P_n}\Big[\big|\widehat \lambda_{n,\eps}\big|^p\big\|\widehat h_{n,\varepsilon}-h_\varepsilon\big\|_{L_{p,\varepsilon}}^p\Big] \\&=:2^{p-1}(I_1+I_2).
  \end{align*}
 The term $I_1$ is controlled by means of Theorem \ref{teolambda} combined with the fact that if $f_{\eps}\in \mathscr F(s,p,q,\mathfrak M_{\eps},A(\eps))$ then $h_{\eps}\in \mathscr F(s,p,q,{\mathfrak M},A(\eps))$, which implies $\|h_{\eps}\|_{L_{p,\eps}}\leq {\mathfrak M}$. 
Concerning the term $I_2$, the Cauchy-Schwarz inequality gives
 \begin{align*}I_2 &\leq \sqrt{\E_{P_n}\big[|\widehat \lambda_{n,\eps}|^{2p}\big]} \int_{A(\eps)}\sqrt{\E_{P_{n}}\big[\big|\widehat h_{n,\varepsilon}(x)-h_\varepsilon(x)\big|^{2p}\big]}dx=\sqrt{J_1}\sqrt{J_2}.\end{align*}
The term $J_1$ is treated using the triangle inequality $\E[|\widehat \lambda_{n,\eps}|^{2p}]\leq C_{p}( \lambda_{\eps}^{2p}+\E[|\widehat \lambda_{n,\eps}-\lambda_{\eps}|^{2p}])$ and Theorem \ref{teolambda}.
For $J_2$, notice that as $A(\eps)$ is bounded, an application of the Jensen inequality yields:
\begin{align*}
\int_{A(\eps)}\sqrt{\E_{P_{n}}\big[\big|\widehat h_{n,\varepsilon}(x)-h_\varepsilon(x)\big|^{2p}\big]}dx\leq C\sqrt{\E_{P_{n}}\big[\|\hat h_{n,\eps}-h\|_{2p}^{2p}\big]},
\end{align*} where $C$ depends on $\overline A$. The rate of the right hand side of the inequality has been studied in Corollary \ref{cor:hB}. \hfill$\Box$

\paragraph*{Acknowledgements} The work of E. Mariucci has been partially funded by the Federal Ministry for Education and Research through the Sponsorship provided by the Alexander von Humboldt Foundation, by
the Deutsche Forschungsgemeinschaft (DFG, German Research Foundation) – 314838170,
GRK 2297 MathCoRe, and by Deutsche Forschungsgemeinschaft (DFG) through grant CRC 1294 'Data Assimilation'.

\bibliographystyle{abbrv}
\bibliography{refs}


\appendix
\section{Additional proofs\label{sec:app}}
\paragraph*{Proof of Lemma \ref{lemma:momenti}}

Using the decompositions \eqref{eq:Xgen} and \eqref{eq:Xito}, for any $\eps\in(0,1]$ it holds:
\begin{align}
\PP(|X_t|\geq \eps)-t\lambda_\eps&=\sum_{n=0}^\infty \PP\Big(\Big|t b_\eps+\Sigma W_t +M_t(\eps)+\sum_{j=1}^nY_j(\eps)\Big|\geq \eps\Big)e^{-\lambda_\eps t}\frac{(\lambda_\eps t)^n}{n!}-t\lambda_\eps\nonumber \\
   &\leq  \PP(|tb_\eps+\Sigma W_t +M_t(\eps)|\geq \eps)+\PP(N_t(\eps)\geq 2)\nonumber\\
   &\quad +\lambda_\eps t  \PP(|tb_\eps+\Sigma W_t +M_t(\eps)+Y_1(\eps)|< \eps)\label{eq:declem1}\\
   &:=T_1+T_2+T_3,\nonumber
\end{align}
for some drift $b_\eps$ that might depend on $\eps$.
We shall complete the proof by showing that $\lim_{t\to 0} \frac{T_j}{t}=0$ for $j=1,2,3$.

 \emph{Control of $T_1$.} By the Markov inequality, for all $k\geq 1$ it holds
$$\PP(|t b_\eps+\Sigma W_t +M_t(\eps)|\geq \eps)\leq \eps^{-2k}\E\Big[\big(t b_\eps +\Sigma W_t + M_t(\eps)\big)^{2k}\Big].$$
Observe that $t b_\eps +\Sigma W_t$ is a Gaussian random variable with mean $t b_\eps$ and variance $t\Sigma^2$ independent of $M_t(\eps)$. Hence, by the binomial theorem, 
\begin{align*}
\E\big[\big(t b_\eps +\Sigma W_t + M_t(\eps)\big)^{2k}\big]&=\E\big[M_t(\eps)^{2k}\big]+\sum_{j=0}^{2k-1} \binom{2k}{j} \E\big[(tb_\eps +\Sigma W_t )^{2k-j}\big]\E\big[M_t(\eps)^j\big]\\
&:=T_{1,1}+T_{1,2}.
\end{align*}
In the sequel, we establish that
$$\lim_{t\to 0} \frac{T_{1,1}}{t}=\mu_{2k}(\eps) \quad \text{and}\quad \lim_{t\to 0}  \frac{T_{1,2}}{t}=0,\quad\forall k\geq1, $$
where $\mu_j(\eps):=\int_{|x|\leq \eps} x^j\nu(dx)$, $j\geq 2$.
To see that, one needs to control the moments of $M_t(\eps)$ and $tb_\eps +\Sigma W_t$.
To compute $\E[M_t(\eps)^j]$ we use that the $k$th cumulant of $M_t(\eps)$ is $t\mu_k(\eps)$ and the properties of Bell's polynomials to derive that for $j\ge2$,
\begin{align}
\E[M_t(\eps)^j]&=B_j(0,t\sigma^2(\eps),t\mu_3(\eps),\dots,t\mu_j(\eps))
= \sum_{l=1}^j t^l B_{j,l}(0,\sigma^2(\eps),\mu_3(\eps),\dots,\mu_j(\eps))\nonumber\\
&=tB_{j,1}(0,\sigma^2(\eps),\mu_3(\eps),\dots,\mu_j(\eps))+O(t^{2})=t\mu_{j}(\eps)+O(t^{2}) \ \text{as }t\to 0, \label{bell}
\end{align}
where $B_j$ and $B_{j,l}$ denote the $j$-th complete exponential Bell polynomial and the incomplete exponential Bell polynomials, respectively.
From \eqref{bell}, it directly follows that $\lim_{t\to 0}\frac{\E[M_t(\eps)]^j}{t}=\mu_j(\eps)$ for all $j\geq 2$ and $\eps\in(0,1]$. Secondly, using the formulas for the moments of Gaussian distributions, we derive that $\lim_{t\to 0} \frac{T_{1,2}}{t}=0$ for all $k\geq 2$ and $\eps>0$. 
Collecting all pieces together, we deduce that
$$\limsup_{t\to 0}\frac{\PP(|t b_\eps+\Sigma W_t +M_t(\eps)|\geq \eps)}{t}\le \frac{\mu_{2k}(\eps)}{\eps^{2k}},\quad \forall k\geq 1.$$

We then conclude that $\lim_{t\to0}\frac{T_1}{t}=0$ by observing that $\lim_{k\to\infty}\frac{\mu_{2k}(\eps)}{\eps^{2k}}=0.$
Indeed, let us write
\begin{align*}
\frac{\int_{-\eps}^{\eps}x^{2k}\nu(dx)}{\eps^{2k}}&=\int_{-\eps}^{\eps}\Big(\frac{x}{\eps}\Big)^{2k}f(x)dx=\Bigg(\int_{-\eps}^{-\eps+\frac{1}{\sqrt k}}+\int_{-\eps+\frac{1}{\sqrt k}}^{\eps-\frac{1}{\sqrt{k}}}+\int_{\eps-\frac{1}{\sqrt k}}^{\eps}\Bigg)\Big(\frac{x}{\eps}\Big)^{2k}f(x)dx\\&=A_{1}(k)+A_{2}(k)+A_{3}(k).
\end{align*}
Using that $\int_{-\eps}^{\eps}x^{2}f(x)dx<\infty$ joined with $\big(\frac{x}{\eps}\big)^{2k}\leq\big( \frac{x}{\eps}\big)^2$, $\forall k\geq 1$ and $x\in[-\eps,\eps]$, we get 
$$A_{1}(k)\leq \frac{1}{\eps^{2}}\int_{-\eps}^{-\eps+\frac{1}{\sqrt k}}x^{2}f(x)dx\to0,\quad \mbox{as }k\to\infty.$$
Similarly, $A_{3}(k)$ tends to 0 as $k\to\infty$. Finally, we have that
\begin{align*}
A_{2}(k)&\leq \Big(1-\frac{1}{\eps\sqrt{k}}\Big)^{2k-2}\frac{\sigma^{2}(\eps)}{\eps^{2}}\sim \frac{\sigma^{2}(\eps)}{\eps^{2}}\exp\Big(-\frac{2k-2}{\eps\sqrt{k}}\Big)\to0,\quad \mbox{as }k\to\infty.
\end{align*} 

\emph{ Control of $T_2$.} As $\PP(N_t(\eps)\geq 2)=1-e^{-\lambda_\eps t}-e^{-\lambda_\eps t}\lambda_\eps t$ it holds  $\lim_{t\to0}\frac{T_2}{t}=0$.

\emph{ Control of $T_3$.} The fact that $\lim_{t\to0}\frac{T_3}{t}=0$ is implied by 
$$\lim_{t\to0}\PP(|tb_\eps+\Sigma W_t +M_t(\eps)+Y_1(\eps)|< \eps)=0$$
which holds true by the dominated convergence theorem as $W_0=0$, $M_0(\eps)=0$ and $|Y_{1}(\eps)|\ge\eps$ a.s.
\hfill$\Box$

\paragraph*{Proof of Lemma \ref{lem:hpeps}}
Using the definition \eqref{eq:peps} we derive that
\begin{align*}
\peps -h_{\eps}&= h_{\eps}\Big(\frac{e^{-\lambda_{\eps}\Delta}\lambda_{\eps}\Delta}{1-e^{-\lambda_{\eps}\Delta}}-1\Big)+\sum_{k=2}^{\infty}\frac{e^{-\lambda_{\eps}\Delta}(\lambda_{\eps}\Delta)^{k}}{k!(1-e^{-\lambda_{\eps}\Delta})}h_{\eps}^{\star k}\\
&=\Big(\frac{e^{-\lambda_{\eps}\Delta}\lambda_{\eps}\Delta}{1-e^{-\lambda_{\eps}\Delta}}-1\Big)h_{\eps}+(\lambda_{\eps}\Delta)^{2}\sum_{k=2}^{\infty}\frac{e^{-\lambda_{\eps}\Delta}(\lambda_{\eps}\Delta)^{k-2}}{k!(1-e^{-\lambda_{\eps}\Delta})}h_{\eps}^{\star k}.
\end{align*} Taking the $L_{p}$ norm and using the Young inequality together with the fact that $h_{\eps}$ is a density with respect to the Lebesgue measure, i.e. $\| h_{\eps}\|_{L_{1,\eps}}\leq 1$, we get \begin{align*}
\big\|\peps -h_{\eps}\big\|_{L_{p,\eps}}&\leq\Big(\Big|\frac{e^{-\lambda_{\eps}\Delta}\lambda_{\eps}\Delta}{1-e^{-\lambda_{\eps}\Delta}}-1\Big|+\frac{(\lambda_{\eps}\Delta)^{2}}{1-e^{-\lambda_{\eps}\Delta}}\Big)\|h_{\eps}\|_{L_{p,\eps}}\leq 2\lambda_\eps\Delta e^{\lambda_\eps\Delta}\|h_{\eps}\|_{L_{p,\eps}},
\end{align*} 
as desired.

\paragraph*{Proof of Lemma \ref{lem:nB}}
We have $\bm n(\varepsilon)=\sum_{i=1}^{n}\1_{(\varepsilon,\infty)}(|X_{i\Delta}-X_{(i-1)\Delta}|)=\widehat \lambda_{n,\eps} n\Delta.$ We introduce the centered i.i.d. random variables $V_i=\1_{(\varepsilon,\infty)}(|X_{i\Delta}-X_{(i-1)\Delta}|)-F_{\Delta}(\eps)$, which are bounded by 2 and such that $\E [V_i^2]\leq F_{\Delta}(\eps)$. Applying the Bernstein inequality we have, 
\begin{align}\label{eq N control}
\PP\Big(\Big|\frac{\nB}{n}-F_{\Delta}(\eps)\Big|>x\Big)
&\leq2\exp\Big(- \frac{nx^2}{2( F_{\Delta}(\eps)+\frac{2x}{3})}\Big),\quad x>0.
\end{align} Fix $x=F_{\Delta}(\eps)/2$, on the set $A_x=\big\{\big|\tfrac{\nB}{n}-F_{\Delta}(\eps)\big|\leq x\big\}$ we have \begin{align}\label{eq prf control1}n\frac{F_{\Delta}(\eps)}{2}\leq \nB\leq n\frac{3F_{\Delta}(\eps)}{2}.\end{align} 
Moreover it holds that
$
\E\big[\nB^{-r}\big]=\E\big[\nB^{-r} \mathds{1}_{{A_x^{c}}}\big]+\E\big[\nB^{-r}\mathds{1}_{A_x}\big].
$ Since $r\geq0$ and $\nB\geq1$, using \eqref{eq N control} and \eqref{eq prf control1} we get the following upper bound
\begin{align*}
\E\big[&\nB^{-r}\big]
\leq 2\exp\big(-\tfrac{3}{32}nF_{\Delta}(\eps)\big)+\Big(\frac{nF_{\Delta}(\eps)}{2}\Big)^{-r}\end{align*} and the lower bound
$
\E\big[\nB^{-r}\big]\geq \E\big[\nB^{-r}\mathds{1}_{A_{x}}\big]\geq \Big(\frac{3nF_{\Delta}(\eps)}{2}\Big)^{-r}.$ 
This completes the proof.\hfill$\Box$

\paragraph*{Proof of Lemma \ref{lem:tnB}}
For the first inequality, the proof is similar to the proof of Lemma \ref{lem:nB}. Using the definition of $\tnB$ we have
$
\tnB=\sum_{i\in\mathscr I_{\eps}}\1_{Z_{i\Delta}(\varepsilon)\ne Z_{(i-1)\Delta}(\varepsilon)}.
$ For $i\in \mathscr I_{\eps}$, we set $W_{i}:=\1_{Z_{i\Delta}(\varepsilon)\ne Z_{(i-1)\Delta}(\varepsilon)}$. We have $$\E[W_{i}|i\in\mathscr I_{\eps}]= \PP\big(Z_{i\Delta}(\varepsilon)\ne Z_{(i-1)\Delta}(\varepsilon)\big||X_{i\Delta}-X_{(i-1)\Delta}|>\eps\big)=1-\frac{v_{\Delta}(\eps)e^{-\lambda_{\eps}\Delta}}{F_{\Delta}(\eps)},$$ using the independence of $M(\varepsilon)$ and $Z(\varepsilon)$. The variables $W_{i}-\E[W_{i}|i\in\mathscr I_{\eps}]$ are centered, i.i.d., bounded by 2 and such that the following bound on the variance holds: $\V(W_{i}|\mathscr I_{\eps})\leq \frac{v_{\Delta}(\eps)e^{-\lambda_{\eps}\Delta}}{F_{\Delta}(\eps)}$.  Applying the Bernstein inequality we have, 
\begin{align}\label{eq:control}
\PP\bigg(\Big|\frac{\tnB}{\nB}-\Big(1-\frac{v_{\Delta}(\eps)e^{-\lambda_{\eps}\Delta}}{F_{\Delta}(\eps)}\Big)\Big|>x\big|\mathscr I_{\eps}\bigg)
&\leq2\exp\bigg(- \frac{\nB x^2}{2\big( \frac{v_{\Delta}(\eps)e^{-\lambda_{\eps}\Delta}}{F_{\Delta}(\eps)}+\frac{2x}{3})\big)}\bigg),\quad x>0.
\end{align} Fix $x=\tfrac{1}{2}$, on the set $A_x=\big\{\big|\tfrac{\tnB}{\nB}-(1-\tfrac{v_{\Delta}(\eps)e^{-\lambda_{\eps}\Delta}}{F_{\Delta}(\eps)})\big|\leq \frac12\big\}$ we have 
\begin{align}
\label{eq control1}\frac\nB6<\Big(\frac12-\frac{v_{\Delta}(\eps) e^{-\lambda_{\eps}\Delta}}{F_{\Delta}(\eps)}\Big)\nB&\leq \tnB,
\end{align} if $\frac{v_{\Delta}(\eps)}{F_{\Delta}(\eps)}\leq \frac13$. It follows from \eqref{eq:control}, \eqref{eq control1} and  $\tnB\geq 1$ that for $r\geq 0$
\begin{align*}
\E\big[\tnB^{-r}\big|\mathscr I_{\eps}\big]&\leq 2\exp\big(-\tfrac{3}{16}\nB \big)+\Big(\frac\nB6\Big)^{-r}.
\end{align*}  Finally, using that for all $x>0$ we have $x^{r}e^{- x}\leq C_{r}:={r}^{r}e^{-r}$, we derive 
\begin{align*}
\E\big[\tnB^{-r}\big|\mathscr I_{\eps}\big]&\leq C_{r} \nB^{-r}+\Big(\frac\nB6\Big)^{-r},\end{align*} which leads to the first part of the result.

The second part of the result can be obtained by means of the Rosenthal inequality. For $r\geq 0$, we have, using that $\nB\geq\tnB$, {\small
\begin{align*}
\E\big[(\nB-\tnB)^{r}\big|\mathscr I_{\eps}\big]&\leq C_{r}\Big\{\E\Big(\Big|\nB\big(1-\tfrac{v_{\Delta}(\eps)e^{-\lambda_{\eps}\Delta}}{F_{\Delta}(\eps)}\big)-\tnB\Big|^{r}\Big|\mathscr I_{\eps}\Big)+\Big(\nB\tfrac{v_{\Delta}(\eps)e^{-\lambda_{\eps}\Delta}}{F_{\Delta}(\eps)}\Big)^{r}\Big\}.
\end{align*}} The Rosenthal inequality leads to, for $r\geq 2$,
\begin{align*}
E\Big[\Big|\nB\big(1-\tfrac{v_{\Delta}(\eps)e^{-\lambda_{\eps}\Delta}}{F_{\Delta}(\eps)}\big)-\tnB\Big|^{r}\Big|\mathscr I_{\eps}\Big]&\leq C_{r}\Big(\nB \tfrac{v_{\Delta}(\eps)e^{-\lambda_{\eps}\Delta}}{F_{\Delta}(\eps)} \Big)^{r/2}.
\end{align*} 
Thanks to the Jensen inequality we can also treat the case $0<r<2$ recovering the same inequality. Therefore, it follows that for all $r>0$
\begin{align*}
\E\big[(\nB-\tnB)^{r}\big|\mathscr I_{\eps}\big]&\leq C\Big\{\Big(\nB\frac{v_{\Delta}(\eps)e^{-\lambda_{\eps}\Delta}}{ F_{\Delta}(\eps)}\Big)^{r/2}+\Big(\nB\frac{v_{\Delta}(\eps)e^{-\lambda_{\eps}\Delta}}{F_{\Delta}(\eps)}\Big)^{r}\Big\}.
\end{align*}  
This completes the proof.
\hfill$\Box$

\paragraph{Proof of Theorem \ref{thm:mainF}}
Fix $0<\eps\le 1$, the proof is a consequence of Theorem \ref{thm:main}. In the sequel $C$ is a constant, possibly depending on $\eps $, $s$, $p$, $\|\Phi\|_{\infty}$, $\|\Phi'\|_{\infty}$, $\|\Phi\|_{p}$ and $\mathfrak M$, whose value may change from line to line. First, using \eqref{eq:assH1} and since $\Delta\le\frac{1}{2c}\wedge\delta\le \frac{\lambda_{\eps}}{2c}$ it holds that $F_{\Delta}(\eps)\ge \lambda_{\eps}\Delta/2$, therefore  $1-F_{\Delta}(\eps)\le e^{-F_{\Delta}(\eps)}\le e^{-\lambda_{\eps}\Delta/2}$. This together with $nv_{\Delta}(\eps)\le 1$, leads to $\frac{v_{\Delta}(\eps)}{F_{\Delta}(\eps)}\le\frac13$ whenever $ n\Delta\ge 6,$ using that $\lambda_{\eps}\ge 1$ on $\mathscr F_{\mathscr H_{1}}$.

Replacing $2^{J}=(n\Delta)^{{1}/{(2s+1)}}$ and using that $n\Delta^{2}\le 1$ and $v_{\Delta}(\eps)\le 1/n$, the upper bound given in Theorem \ref{thm:main} can be rewritten in 
 \begin{align*}
  \big[\ell_{p,\varepsilon}\big(\widehat f_{n,\varepsilon},f\big)\big]^p&\leq C\bigg\{ e^{-n\Delta\lambda_{\eps}/2}+\Big(\frac{1}{n\Delta}\Big)^{\frac{p}{2}}+\Delta^p
+\frac{(n\Delta)^{\frac{2p}{2s+1}}}{(n \Delta)^{p}}+(n\Delta)^{-\frac{sp}{2s+1}}\\
&\hspace{-1cm}+(n\Delta)^{\frac{-2s(p-1)}{2s+1}}\mathbf{1}_{p\ge 2}
 +(n\Delta)^{\frac{5p-2}{2(2s+1)}}\Big[(n{\Delta})^{1-p} \big(\Delta+\Delta^{\frac p2}\big)+n^{-\frac p2}+\Delta^p\Big]\bigg\}\\
 &\hspace{-1cm}\le C\bigg\{ 
(n\Delta)^{-\frac{sp}{2s+1}}+(n\Delta)^{\frac{p(1-2s)}{(2s+1)}}
+(n\Delta)^{\frac{5p-2}{2(2s+1)}}\Big[(n{\Delta})^{1-p} \big(\Delta+\Delta^{\frac p2}\big)\Big]\bigg\},
  \end{align*} 
where as $\eps$ is fixed, the quantities $\lambda_{\eps}$, $\ \mu_{p}(\eps),\ \sigma^{2}(\eps)$ and $b_\nu(\eps)$ are included in the constant $C$, together with the terms $\|f\|_{p/2,\eps}\le \|f\|_{p,\eps}^{2}\le (\lambda_{\eps}\mathfrak M)^{2}$ using the H\"older inequality and that $f\in \mathscr	F({s,p,q,\mathfrak{M}_\eps,A(\eps)})$. For $p\ge 2$, i.e. $\Delta^{p/2}\le\Delta$, the dominating terms in the latter inequalities are the following
  \[v^{*}:=(n\Delta)^{-\frac{sp}{2s+1}},\quad v_{1}:=(n\Delta)^{\frac{p(1-2s)}{2s+1}},\quad v_{2}:=(n\Delta)^{\frac{3p-4s(p-1)}{2(2s+1)}}\Delta.  \]
The following computations lead, for $s\ge \frac32-\frac1p$ and $p\ge2$, to Theorem \ref{thm:mainF} as
  \begin{align*}
   \frac{v_{1}}{v^{*}}&=(n\Delta)^{\frac{p(1-s)}{2s+1}}\le 1,\quad
\frac{v_{2}}{v^{*}}=(n\Delta)^{\frac{3p-2s(p+1)}{2(2s+1)}}n\Delta^{2}\le 1\quad
  \mbox{and}\quad \frac{v_{2}}{v_{1}}&=(n\Delta)^{\frac{p+4s}{2(2s+1)}}\Delta .&
  \end{align*}

  \paragraph{Proof of Theorem \ref{thm:mainF2}}
Fix $0<\eps\le 1$, the proof is a consequence of Theorem \ref{thm:main}. In the sequel $C$ is a constant, possibly depending on $\eps $, $s$, $p$, $\|\Phi\|_{\infty}$, $\|\Phi'\|_{\infty}$, $\|\Phi\|_{p}$ and $\mathfrak M$, whose value may change from line to line. First, using \eqref{eq:assH1} and that $\lambda_{\eps}\ge 1$, it holds for $\Delta\le\frac{1}{6c}\wedge\delta\le \frac{\lambda_{\eps}}{2c}$ that $F_{\Delta}(\eps)\ge \lambda_{\eps}\Delta/2$ and $1-F_{\Delta}(\eps)\le e^{-\lambda_{\eps}\Delta/2}$. Additionally, using the latter together with \eqref{eq:assH2}, $\beta>1$, and $\Delta\le \frac{1}{6c}$ gives $\frac{v_{\Delta}(\eps)}{F_{\Delta}(\eps)}\le\frac13$.
Replacing $2^{J}=(n\Delta)^{{1}/{(2s+1)}}$, the upper bound given in Theorem \ref{thm:main} can be rewritten 
 \begin{align}
  \big[\ell_{p,\varepsilon}\big(\widehat f_{n,\varepsilon},f\big)\big]^p&\le C\bigg\{  e^{-n\Delta\lambda_{\eps}/2}+\Big(\frac{1}{n\Delta}\Big)^{\frac{p}{2}}+\Delta^p
 +(n\Delta)^{\frac{2p}{2s+1}}\Big[\Big(\frac{\Delta^{\beta-1}}{n \Delta}\Big)^{p/2}+\Delta^{p(\beta-1)}\Big]\nonumber\\
 & \hspace{-2.3cm}+(n\Delta)^{-\frac{sp}{2s+1}}+(n\Delta)^{-\frac{2s(p-1)}{2s+1}}\mathbf{1}_{p\ge 2}
 +(n\Delta)^{\frac{5p-2}{2(2s+1)}}\Big[(n\Delta)^{1-p} \big(\Delta+\Delta^{\frac p2}\big)+n^{-\frac p2}+\Delta^p\Big]\bigg\}\nonumber\\
  &\leq C\bigg\{(n\Delta)^{-\frac{sp}{2s+1}} +(n\Delta)^{\frac{2p}{2s+1}}\Big[\Big(\frac{\Delta^{\beta-1}}{n \Delta}\Big)^{p/2}+\Delta^{p(\beta-1)}\Big]\nonumber\\
 &\quad\quad+(n\Delta)^{\frac{5p-2}{2(2s+1)}}\Big[(n\Delta)^{1-p} \big(\Delta+\Delta^{\frac p2}\big)+n^{-\frac p2}+\Delta^p\Big]\bigg\},\label{eq:prf331}\end{align}
where, $\eps$ being fixed, the quantities  $\lambda_{\eps}, $ $\ \mu_{p}(\eps),\ \sigma^{2}(\eps)$ and $b_\nu(\eps)$ are included in the constant $C$ as well as  $\|f\|_{p/2,\eps}\le \|f\|_{p,\eps}^{2}\le (\lambda_{\eps}\mathfrak M)^{2}$ using the H\"older inequality and that $f\in \mathscr	F({s,p,q,\mathfrak{M}_\eps,A(\eps)})$.
Consider the case $n\Delta^{2}> 1$,   Equation \eqref{eq:prf331} simplifies in 
\begin{align*}
 \big[\ell_{p,\varepsilon}\big(\widehat f_{n,\varepsilon},f\big)\big]^p&\le C\bigg\{(n\Delta)^{-\frac{sp}{2s+1}} +(n\Delta)^{\frac{2p}{2s+1}}\Big[\Big(\frac{\Delta^{\beta-1}}{n \Delta}\Big)^{p/2}+\Delta^{p(\beta-1)}\Big]+(n\Delta)^{\frac{5p-2}{2(2s+1)}}\Delta^p\bigg\}.
\end{align*}
 Set $v^{*}:=(n\Delta)^{-\frac{sp}{2s+1}}$ and
  \[  v_{1}:=(n\Delta)^{\frac{3p-2sp}{2(2s+1)}}\Delta^{\frac p2(\beta-1)},\ v_{2}:=(n\Delta)^{\frac{2p}{2s+1}}\Delta^{p(\beta-1)},\ v_{3}:=(n\Delta)^{\frac{5p-2}{2(2s+1)}}\Delta^p.\]
Next, note that $  \frac{v_{1}}{v_{2}}=(n\Delta^{\beta})^{-\frac p2}$ and \begin{align*}
    \frac{v_{1}}{v^{*}}&=(n\Delta^{\beta})^{\frac{3p}{2(2s+1)}}\Delta^{\frac{p(\beta-1)}{2s+1}(s-1) }=(n\Delta^{2})^{\frac{3p}{2(2s+1)}}\Delta^{\frac{p}{2(2s+1)}(2s(\beta-1)+\beta-4) },\\
   \frac{v_{2}}{v^{*}}&=(n\Delta^{\beta})^{\frac{2p+sp}{2s+1}}\Delta^{ps\frac{(\beta-1)}{2s+1}}  =(n\Delta^{2})^{\frac{2p+sp}{2s+1}}\Delta^{p\beta-3p\frac{1+s}{2s+1}}, \\  \frac{v_{3}}{v^{*}}&=(n\Delta^{2})^{\frac{5p+2sp-2}{2(2s+1)}}\Delta^{\frac{2sp-3p+2}{2(2s+1)}}\ge 1\hspace{2cm} \mbox{if }s\le\frac32-\frac1p,\\
   \frac{v_{2}}{v_{3}}&=(n\Delta)^{\frac{2-p}{2(2s+1)}}\Delta^{p(\beta-2)}\quad \mbox{and}\quad\frac{v_{1}}{v_{3}}=(n\Delta^{2})^{\frac{1-p-sp}{2s+1}}\Delta^{\frac p2(\beta-2)+\frac{p-2}{2(2s+1)}}.
 \end{align*} It follows that, if $\beta>2$ and as $p\ge 2$ then $v_{3}\ge v_{2}$ and $v_{3}\ge v_{1}$. 
If $\beta\in(1,2)$, the constraint $n\Delta^{2}>1$ implies $n\Delta^{\beta}\ge 1$ and ${v_{1}}<v_{2}$. The order of the other terms depends on the rate of $\Delta$ according to $n$.

For the case $n\Delta^{2}\le 1$,  the case $\beta\ge2$ is covered by Theorem \ref{thm:mainF}. If $\beta\in(1,2)$, we add the terms $v_{1}$ and $v_{2}$ that now intervene in the rate. Theorem \ref{thm:mainF2} follows.
  
    \paragraph{Proof of Theorem \ref{thm:mainF3}} 
  
Consider $f\in \mathscr	 L_{M,\alpha}$ for $\alpha\in(0,2)$ and $M>0$; straightforward computations give for $r\in\{1,p/2,p\}${\small,
\begin{align}\label{eq:maj} \int_{|x|>\eps}f(x)^r dx\le \frac{2M^{r}}{r-1+r\alpha}\eps^{1-r-r\alpha} + \int_{|x|>1}f(x)^r dx,\quad \mu_{p}(\eps)\le\frac{2M}{p-\alpha}\eps^{p-\alpha}+\mu_{p}(1).\end{align}}
  
  For $\alpha\in(0,1)$, set $\mathbf{C}_{1,\alpha}:=\max\{2{\rm C}_{1}, {\rm D}_{1}+{\rm D}_{1}(2/\alpha+\lambda_{1})+2(2/\alpha+\lambda_{1})^{2}\}$ where ${\rm C}_{1},\ {\rm D}_{1}$ and ${\rm D}_{2}$ appear in Theorems \ref{thm:vd1} and \ref{thm:Fd1}. Theorem \ref{thm:Fd1} and \eqref{eq:maj} give for $\Delta<\frac{(1-\alpha)\eps^{\alpha}}{M4^{1+\alpha}}$ that $|F_{\Delta}(\eps)-\lambda_{\eps}\Delta|\le M^{2}{\mathbf{C}}_{1,\alpha}\Delta^{2}\eps^{-2\alpha}$, using that $\lambda_{\eps}\ge 1$ and $\Delta\le \frac{\eps^{2\alpha}}{2M^{2}{\mathbf{C}}_{1,\alpha}}$ it holds $\frac{\lambda_{\eps}\Delta}{2}\le F_{\Delta}(\eps)$ and $(1-F_{\Delta}(\eps))\le e^{-\lambda_{\eps}\Delta/2}.$
    Additionally, it follows from Theorem \ref{thm:vd1} that  $v_{\Delta}(\eps)\le M^{2}\mathbf{C}_{1,\alpha}\Delta^{2}\eps^{-2\alpha}$. Thus, for $\Delta \le  \frac{\eps^{2\alpha}}{6M^{2}{\mathbf{C}}_{1,\alpha}}$ it follows that $\frac{v_{\Delta}(\eps)}{F_{\Delta}(\eps)}\le \frac13$. 
    
    For $\alpha\in[1,2)$, set $\mathbf{C}_{2,\alpha}:={\rm F}_{1}+2\alpha{\rm F}_{2}+2{\rm F}_{3}(2/\alpha+\lambda_{1})^{2}$ where ${\rm F}_{1},\ {\rm F}_{2}$ and ${\rm F}_{3}$ appear in Theorem  \ref{thm:Fd2}. Theorem \ref{thm:Fd2} and \eqref{eq:maj} give for $\Delta< \frac{(2-\alpha)\eps^\alpha}{2^{1+\alpha}M}$ that $|F_{\Delta}(\eps)-\lambda_{\eps}\Delta|\le M^{2}{\mathbf{C}}_{2,\alpha}\Delta^{2}\eps^{-2\alpha}$. Similarly, using that $\lambda_{\eps}\ge 1$ and $\Delta\le \frac{\eps^{2\alpha}}{2M^{2}{\mathbf{C}}_{2,\alpha}}$ it holds $\frac{\lambda_{\eps}\Delta}{2}\le F_{\Delta}(\eps)$ and $(1-F_{\Delta}(\eps))\le e^{-\lambda_{\eps}\Delta/2}.$
    Additionally, it follows from \eqref{eq:Mplus} that  $v_{\Delta}(\eps)\le M^{2}\mathbf{C}_{2,\alpha}\Delta^{2}\eps^{-2\alpha}$. Hence, for $\Delta \le  \frac{\eps^{2\alpha}}{6M^{2}{\mathbf{C}}_{2,\alpha}}$ we deduce that $\frac{v_{\Delta}(\eps)}{F_{\Delta}(\eps)}\le \frac13$. 

    We set
    \begin{align}
    \label{eq:Delta34}\overline{\Delta}:=\begin{cases}
    1\wedge \frac{(1-\alpha)\eps^{\alpha}}{M4^{1+\alpha}} \wedge \frac{\eps^{2\alpha}}{6M^{2}{\mathbf{C}}_{1,\alpha}} &\mbox{if }\alpha\in(0,1),\\
   1\wedge  \frac{(2-\alpha)\eps^\alpha}{2^{1+\alpha}M}\wedge  \frac{\eps^{2\alpha}}{6M^{2}{\mathbf{C}}_{2,\alpha}}            &\mbox{if }\alpha\in[1,2).
    \end{cases}
    \end{align}
    In the sequel $C$ is a constant, possibly depending on $\alpha$, $M$, $s$, $p$, $\|f\|_{L_{p,1}}$, $\|f\|_{L_{p/2,1}},$ $\lambda_{1},$  $\mu_{p}(1)$, $\|\Phi\|_{\infty}$, $\|\Phi'\|_{\infty}$, $\|\Phi\|_{p}$ and $\mathfrak M$, whose value may change from line to line. 

     The proof is a consequence of Theorem \ref{thm:main}. First, we choose the resolution $J$ that perform the compromise between the terms
 $2^{-Jsp}$ and $2^{Jp/2}\ell_{p,\eps}^{p/2}\big(n F_{\Delta}(\eps)\big)^{-p/2},$ where  $\lambda_{\eps}\ell_{p,\eps}=\|f_{\eps}\|_{L_{p/2,\eps}}\le C_{p,M,\alpha,\|f\|_{p/2,1}}\eps^{\frac{2}{p}-(1+\alpha)} $ from Equation \eqref{eq:maj}. 
 We choose $J$ such that $2^{J}=(\eps^{1+\alpha-\frac2p}n\Delta\lambda_{\eps}^{2})^{{1}/{(2s+1)}}$ and replace it in the bound of Theorem \ref{thm:main} that  can be rewritten for $p\ge 2$ and after simplification in
\begin{align*}
  \big[\ell_{p,\varepsilon}\big(\widehat f_{n,\varepsilon},f\big)\big]^p&\leq C\bigg\{\Big(\frac{\lambda_\eps}{n\Delta}\Big)^{\frac{p}{2}}+\Delta^p\eps^{-2\alpha p}+\lambda_{\eps}^{p}(\eps^{1+\alpha-\frac2p}n\Delta\lambda_{\eps}^{2})^{-\frac{sp}{2s+1}} \\  & \hspace{-1cm} +\lambda_{\eps}^{p}{(\eps^{1+\alpha-\frac2p}n\Delta\lambda_{\eps}^{2})^{\frac{2p}{2s+1}}}\Big[\Big(\frac{\eps^{-2\alpha}}{n \lambda_{\eps}^{2}}\Big)^{p/2}+\Big(\Delta\eps^{-2\alpha}{\lambda_{\eps}^{-1}}\Big)^{p}\Big] +(\Delta\lambda_{\eps})^{p}\eps^{1-p-p\alpha}
\\
 & \hspace{-1cm}+\lambda_{\eps}^{p}{(\eps^{1+\alpha-\frac2p}n\Delta\lambda_{\eps}^{2})^{\frac{5p-2}{2(2s+1)}}}\Big[(n\Delta\lambda_{\eps})^{1-p} \big(\Delta\eps^{p-\alpha}+(\Delta\eps^{2-\alpha})^{\frac p2}\big) \\
 &\hspace{2cm} +(\lambda_{\eps}n\Delta)^{-\frac p2}\big(\eps^{2-\alpha}\Delta\big)^{\frac p2}+(\eps^{1-\alpha}\Delta)^p\Big]\bigg\}\\
 :&= v_{1}+v_{2}+v^{*}+v_{3}+v_{4}+v_{5}+v_{6}+v_{7}+v_{8}+v_{9}.
  \end{align*}
Using the assumptions $\eps^{1+\alpha-\frac2p}n\Delta\lambda_{\eps}^{2}\ge 1$, $n(\Delta\lambda_{\eps}\eps^{-2\alpha})^{2}\le 1$, $\Delta\le 1$, $\eps\le 1$, $\lambda_{\eps}\ge 1$ and $\lambda_{\eps}\eps^{\alpha}\le 1$, we obtain for $p\ge 2$ and $s\ge \frac32-\frac1p$ that
  \begin{align*}\frac{v_{2}}{v_{1}}&=(n(\Delta\lambda_{\eps}\eps^{-2\alpha})^{2})^{p/2}\Delta^{p/2}\lambda_{\eps}^{-3p/2}\le 1,\\
  \frac{v_{1}}{v^{*}}  &=(\eps^{1+\alpha-\frac2p} n\Delta\lambda_{\eps}^{2})^{{p}(\frac{s}{2s+1}-\frac12)}{\lambda_{\eps}}^{\frac p2}\eps^{\frac p2(1+\alpha)-1}\le 1,\quad\quad
 \frac{v_{3}}{v_{4}}={(n\Delta^{2}\eps^{-2\alpha})^{-p/2}}>1, \\
\frac{v_{3}}{v^{*}}&=(\eps^{-4\alpha}n\Delta^{2}\lambda_{\eps}^{2})^{\frac{p(2+s)}{2(2s+1)}}(n\lambda_{\eps}^{2})^{-p/2}\eps^{4\alpha\frac{p(2+s)}{2s+1}+(1+\alpha-\frac2p)\frac{p(2+s)}{2s+1}-\alpha p},\\
&\le (\eps^{-4\alpha}n\Delta^{2}\lambda_{\eps}^{2})^{\frac{p(2+s)}{2(2s+1)}}(n\lambda_{\eps}^{2})^{-p/2}\eps^{5\alpha\frac{p(2+s)}{2s+1}-\alpha p}\le 1,\\
 \frac{v_{5}}{v^{*}}
&=(n\Delta^{2}\lambda_{\eps}^{2}\eps^{-4\alpha})^{\frac{sp}{2s+1}}(\Delta\eps^{-1-\alpha})^{\frac{sp+p}{2s+1}}\eps^{\frac{1+4\alpha{sp}}{2s+1}}\le 1,\\
  \frac{v_{6}}{v_{7}}&=(\Delta\eps^{-\alpha})^{1-\frac p2} \ge 1, \quad\quad\frac{v_{6}}{v_{8}}=(n\Delta^{2}\lambda_{\eps}\eps^{-\alpha})^{1-{\frac p2}}\ge 1,\quad\quad
  \frac{v_{6}}{v_{9}}=(n\Delta^{2}\lambda_{\eps}\eps^{-\alpha} )^{1-p}\ge 1,\\
 \mbox{and }\quad \frac{v_{6}}{v^{*}}&=(\eps^{1+\alpha-\frac2p}n\Delta\lambda_{\eps}^{2})^{\frac{3p-2-2sp}{2(2s+1)}}(n\Delta^{2}\lambda_{\eps}^{2}\eps^{-4\alpha})(\lambda_{\eps}\eps^{2})^{p-1}\eps^{3\alpha+\alpha p}\le 1.
  \end{align*} The result follows.

\section{Some inequalities on the class $\mathscr	L_{M,\alpha}$\label{app:note}}

We partially reproduce here the main results of \cite{DMnote} that provide a control of the quantities $F_{\Delta}(\eps)$ and $v_{\Delta}(\eps)$ on the class $\mathscr	L_{M,\alpha}$, $\alpha\in(0,2),\ M>0$  (see \eqref{eq:Lalpha}). 
Hereafter, the dependency in $\alpha$ of the constants is not given. For explicit values of the constants in the following Theorems, the reader is referred to \cite{DMnote}.

\begin{theorem}\label{thm:vd1}
Let $\nu$ be a Lévy measure absolutely continuous with respect to the Lebesgue measure and denote by $f=\frac{d\nu}{dx}$. Let $\eps\in(0,1]$, $\alpha\in(0,1)$, $M>0$, $f\in\mathscr L_{M,\alpha}$ and $t\in(0,(1-\alpha) M^{-1} \eps^{\alpha}4^{-(1+\alpha)})$. Then, there exists a constant $  {\rm C}_1>0$, only depending on $\alpha$, such that
$\PP(|tb_{\nu}(\eps)+M_t(\eps)|\geq \eps)\leq 2t^2 M^2 {\rm C}_1 \eps^{-2\alpha}.$
\end{theorem}

\begin{theorem}\label{thm:Fd1}
Let $X$ be a finite variation Lévy process of the form $X_t=\sum_{0<s\leq t}\Delta X_s$ with Lévy measure $\nu$ absolutely continuous with respect to the Lebesgue measure and denote by $f=\frac{d\nu}{dx}$. Suppose that $f\in\mathscr L_{M,\alpha}$ for some $\alpha\in(0,1)$ and $M>0$. If $\eps\in(0,1]$, then there exist two constants $\rm D_1$ and $\rm D_2$  only depending on $\alpha$, such that for all $t\in(0,(1-\alpha) M^{-1} \eps^{\alpha}4^{-(1+\alpha)})$ it holds:
\begin{equation*}
|\PP(|X_t|>\eps)-\lambda_\eps t|\leq t^2\big(M^2\eps^{-2\alpha} {\rm D_1} +M\lambda_\eps\eps^{-\alpha} \rm D_2+2\lambda_\eps^2\big).
\end{equation*}
\end{theorem}

\begin{theorem}\label{thm:vd2}
Let $X$ be a Lévy process as in \eqref{eq:defX} and let $\nu$ be a symmetric Lévy measure absolutely continuous with respect to the Lebesgue measure and denote by $f=\frac{d\nu}{dx}$. Let $\eps\in(0,1]$, $\alpha\in[1,2)$, $M>0$, $f\in\mathscr L_{M,\alpha}$ and $t\in(0,(\eps/2)^\alpha(1\land ((2-\alpha)/2M))$. Then, there exists a constant $ {\rm E}_{1}>0$, only depending on $\alpha$, such that
\begin{align*}
\PP(|M_t(\eps)|\geq \eps)&\leq\frac{2^{2+\alpha}M t^{1+1/\alpha}}{\eps^{1+\alpha}}\bigg(1+\frac{M}{\alpha(2-\alpha)(\alpha-1)}\bigg)+ 2t^2 M^2 {\rm E}_{1} \eps^{-2\alpha}, \quad \alpha\in(1,2),\\
\PP(M_t(\eps)\geq \eps)&\leq  \frac{4t^2M^2}{\eps^2}\bigg(e^{2+1/e}+ \frac{37}{9}\bigg)+\frac{4M t^2}{\eps^2}+\frac{16M^2}{\eps^2}t^2 \ln\Big(\frac{\eps}{2t}\Big), \quad \alpha=1.
\end{align*}
\end{theorem}

\begin{theorem}\label{thm:Fd2bis}
Let $X$ be a Lévy process as in \eqref{eq:defX} and let $\nu$ be a symmetric L\'evy measure with density $f$ with respect to the Lebesgue measure and $f\in \mathscr L_{M,\alpha}$ for some $\alpha\in[1,2)$ and $M>0$. Then, for all $0<t<(\eps/2)^\alpha\big(1\land ((2-\alpha)/2M)\big)$, $\eps\in(0,1]$, it holds: 
\begin{align*}
|\PP(|X_t|>\eps)-&\lambda_\eps t|\leq {\rm G}_{1}\frac{t^{1+1/\alpha}}{\eps^{1+\alpha}}+{\rm G}_{2}\frac{t^{2}}{\eps^{2\alpha}}+{\rm G}_3\frac{t^{2}}{\eps^{2}} \ln\Big(\frac{\eps}{t}\Big)\1_{\alpha=1},\end{align*}
where  ${\rm G}_{1}$, ${\rm G}_{2}$  and ${\rm G}_{3}$ are positive constants, only depending on $M$, $\alpha$ and $\lambda_{1}$.

\end{theorem}

\begin{theorem}\label{thm:Fd2}
Let $X$ be a Lévy process as in \eqref{eq:defX} and let $\nu$ be a symmetric L\'evy measure having a density $f$ with respect to the Lebesgue measure with $f\in \mathscr L_{M,\alpha}$ for some $\alpha\in[1,2)$ and $M>0$. Let $\eps\in(0,1]$ and assume that $f$ is  $M\eps^{-(2+\alpha)}$-Lipschitz on the interval $(3/4\eps,5/4\eps)$. For all $t\in (0, \frac{(2-\alpha)\eps^\alpha}{2^{1+\alpha}M})$, it holds:\begin{align*}
|\PP(|X_t|>\eps)-\lambda_\eps t|&\leq
 t^2M^2 \big( {\rm F}_{1}\eps^{-2\alpha}+\lambda_1\eps^{-\alpha}{\rm F}_2 \big)+2t^2\lambda_1^2 +\frac{t^4 M^4 {\rm F}_3}{\eps^{4\alpha}},
\end{align*}
where ${\rm F}_1$, ${\rm F}_2$ and ${\rm F}_3$ are universal positive constants, only depending on $\alpha$.
\end{theorem}
Note that Theorem \ref{thm:Fd2} applied to $\nu=\nu\mathbf{1}_{[-\eps,\eps]}$, satisfying the assumptions of Theorem \ref{thm:Fd2}, permits to improve the result of  Theorem \ref{thm:vd2}  as follows
\begin{align}\label{eq:Mplus}
\PP(|M_t(\eps)|>\eps)&\leq
 M^2  {\rm F}_{1}t^2\eps^{-2\alpha}+\frac{t^4 M^4 {\rm F}_3}{\eps^{4\alpha}}.
\end{align}

If $\alpha\in(0,1)$, Theorem \ref{thm:Fd1} permits to derive that Assumption \eqref{eq:assH1} is valid for $\delta=\frac{(1-\alpha)\eps^{\alpha}}{M4^{1+\alpha}}=:\delta_{0}$ and a constant $c$ depending only on $\alpha$, $M$ and whose dependency in $\eps$ is explicit. Theorem \ref{thm:vd1} ensures that Assumption \eqref{eq:assH2} is fulfilled with $\beta=2$, $\delta=\delta_{0}$ and $c$ depending only on $\alpha$, $M$ and whose dependency in $\eps$ is explicit. 

If $\alpha\in[1,2)$ and the Lévy measure $\nu$ is symmetric, Theorem \ref{thm:Fd2} (under an addition local Lipschitz condition on $f$) permits to derive that Assumption \eqref{eq:assH1} is valid for $\delta=\frac{(2-\alpha)\eps^{\alpha}}{M2^{1+\alpha}}=:\delta_{1}$ and a constant $c$ depending only on $\alpha$, $M$ and whose dependency in $\eps$ is explicit. Theorem \ref{thm:vd2} ensures that Assumption \eqref{eq:assH2} is fulfilled with $\beta=1+1/\alpha$, $\delta=\delta_{1}$ (if $M\ge1/2$). This can be improved using \eqref{eq:Mplus} in $\beta=2$ under a local Lipschitz condition on $f$.

\end{document}